\documentclass[a4paper, 12pt]{article}

\usepackage{amsmath,amssymb,amsfonts,setspace}
\usepackage{mathtools}

\usepackage{authblk}

\newcommand{\eq}{\begin{equation}}
\newcommand{\en}{\end{equation}}

\def\bp{\boldsymbol p}

\def\br{\boldsymbol r}

\newtheorem{theorem}{Theorem}

\newtheorem{lemma}{Lemma}

\newtheorem{corollary}{Corollary}

\newtheorem{prop}{Proposition}

\def\endpf{\hfill $\Box$ \vskip0.5cm}
\def \proof{\noindent{\it Proof.\ }}

\usepackage{adjustbox}

\newcommand{\veq}{\mathrel{\rotatebox{90}{$=$}}}
\newcommand{\vless}{\mathrel{\rotatebox{90}{$<$}}}
\newcommand{\vmore}{\mathrel{\rotatebox{90}{$>$}}}

\newlength{\bibitemsep}
\newlength{\bibparskip}
\let\oldthebibliography\thebibliography
\renewcommand\thebibliography[1]{%
  \oldthebibliography{#1}%
  \setlength{\parskip}{\bibitemsep}%
  \setlength{\itemsep}{\bibparskip}%
}

\begin{document}

\author{Alexander Gnedin }
\affil{Queen Mary, University of London}
\title{Cross Modality of the Extended Binomial Sums}

\date{}

\pagestyle{plain}

\maketitle

\begin{abstract}
\noindent
For a family of probability functions (or a probability kernel), cross modality occurs when every  likelihood maximum matches a mode of the distribution.
This implies existence of simultaneous maxima on the modal ridge of the family.
The paper explores the property for  extended Bernoulli sums, which are  random variables representable as  a sum of  independent  Poisson and any number (finite or infinite) of  Bernoulli
random variables with variable success probabilities.
We show that the cross modality holds for many subfamilies of the class, including  power series distributions derived from 
entire functions with totally positive series expansion.
A central role in the study is played by the extended  Darroch's rule \cite{Darroch,  Samuels}, which originally  
 localised the mode of Poisson-binomial distribution in terms of the mean.
We give different proofs and  geometric interpretation to the extended rule and point at other modal properties of extended Bernoulli sums,
in particular discuss stability of the mode in the context of a transport problem.

\end{abstract}

\section{Itroduction}
A starting point for the  {\it cross modality} phenomenon discussed in this paper is  the following simple observation
on a sequence of   independent identical Bernoulli trials. 
Let $S_n$ be the number of successes in  the first $n$ trials, 
and $T_k$  be the index of  the $k$th success.
Then,
 {for $\ell$   chosen  as the most likely value of $T_{k+1}-1$,  the most likely value of  
$S_{\ell}$ will be $k$}. 
The distribution of $T_{k+1}-1$ is proportionate to the likelihood function of $S_n$ with $n$ regarded as unknown parameter.
Therefore the pair $(k, \ell)$ appears as a maximiser for both the probability function of the number of successes and the likelihood.
In combinatorial terms,  viewing the binomial probabilities as entries of the normalised Pascal triangle, 
we have that every column maximum is also a row maximum.
It is natural to wonder to which extent this kind of property is common.

To put the question in a formal framework and introduce some terminology,
let $f(x;t)$ be a  family of probability functions (continuous or discrete densities, possibly representing a Markov kernel `from $t$ to $x$'), where $x$ runs over the range  $R$ and $t$ is a parameter.
 We call  $(x_0;t_0)$ a {\it cross mode} if the maximum of probability $f(x; t_0)$ on $R$ is attained at $x_0$, and the maximum of likelihood $f(x_0; t)$ is attained at $t_0$.  The family  is said to be {\it cross modal} if, 
for each  $x_0$ in $R$,  every  maximum point $t_0$ of the likelihood function $f(x_0; t)$ is such that $(x_0; t_0)$ is a cross mode. The maxima in these definitions are assumed to exist and meant in the  absolute sense.

We call {\it  peak height} (or just peak) $h(t)$ the modal probability, that is  the maximum value  of  the probability function for fixed parameter $t$, and denote  $m(t)$ the set of {\it modes} where the maximum  is attained. 
In general $m(t)$ need not  be a singleton,  
in particular  univariate unimodal distributions  may have an interval of modes
 \cite{Joag-Dev}.
Using  analogy with  a mountainous landscape, we call
the  `curve'  $M=\{(x,t): x\in m(t)\}$   
 the {\it modal ridge} of the family.
The cross modality means  that each $x$ paired with any of its associated likelihood maximisers belongs to the  modal ridge.
Denoting $\ell(x)$ the set of likelihood maximisers for given $x$, the cross modality condition can be expressed  as an  inversion relation for multivalued functions:
$$ x\in \bigcap_{t\in \ell(x)} m(t).$$
This simplifies for single-valued $m$ as that $t\in\ell(x)$ implies $m(t)=x$, and further simplifies as 
$m(\ell(x))=x$ if also $\ell$ is single-valued.

\begin{figure}
\centering
\includegraphics[scale=0.4]{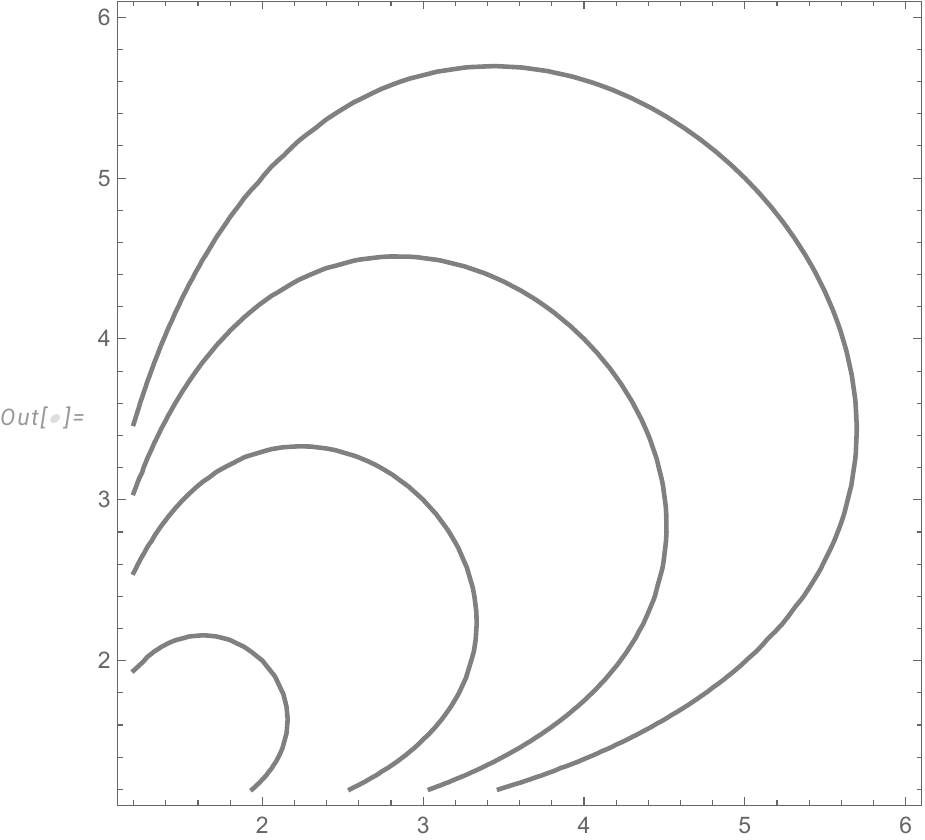}
\caption{Contours of the peak height in the $(a,b)$-plain for the beta distribution $f(x; a,b)= x^{a-1}(1-x)^{b-1}/{\rm B}(a,b), x\in [0,1]; a,b\geq 1$.
 }
\label{BETA}
\end{figure}

The cross modality requires the
 parameter space  to be sufficiently rich in order to match, by means of the likelihood maximisation,  every  $x_0$ in $R$  with a mode
of some distribution from the family. 
Obvious examples are   the location families of the form $f(x;t)=g(x-t)$  (with $x,t \in {\mathbb R}^d$), where $g$ is a unimodal probability function. 
However, 
for general continuous distributions with continuous parameters
the property   is more of an exception than a rule.
In particular, when both $x$ and $t$  are one-dimensional, 
one may expect $M$ to be a planar curve
which has segments where the peaks are of the same height. For,  
at each cross mode
the  gradient  vanishes 
 hence  $f$ has zero derivative in the direction of the modal ridge. With regards to genericity
this degree of  instability under small deformations of the shape
is comparable with  that for  distributions having an interval of modes. 
In the differential topology, such functions having  a manifold of critical points belong to the scope
 of the   Morse-Bott theory \cite{Bott}.
Guided by this analogy, it is possible to achieve the cross modality within  the classic  multiparameter families of distributions (e.g.  gamma, beta or multivariate normal), 
by restricting their parameters to vary along   peculiar contours where the peaks are constant. See Figure \ref{BETA} for illustration.

\begin{figure}
\centering
\includegraphics[scale=0.5]{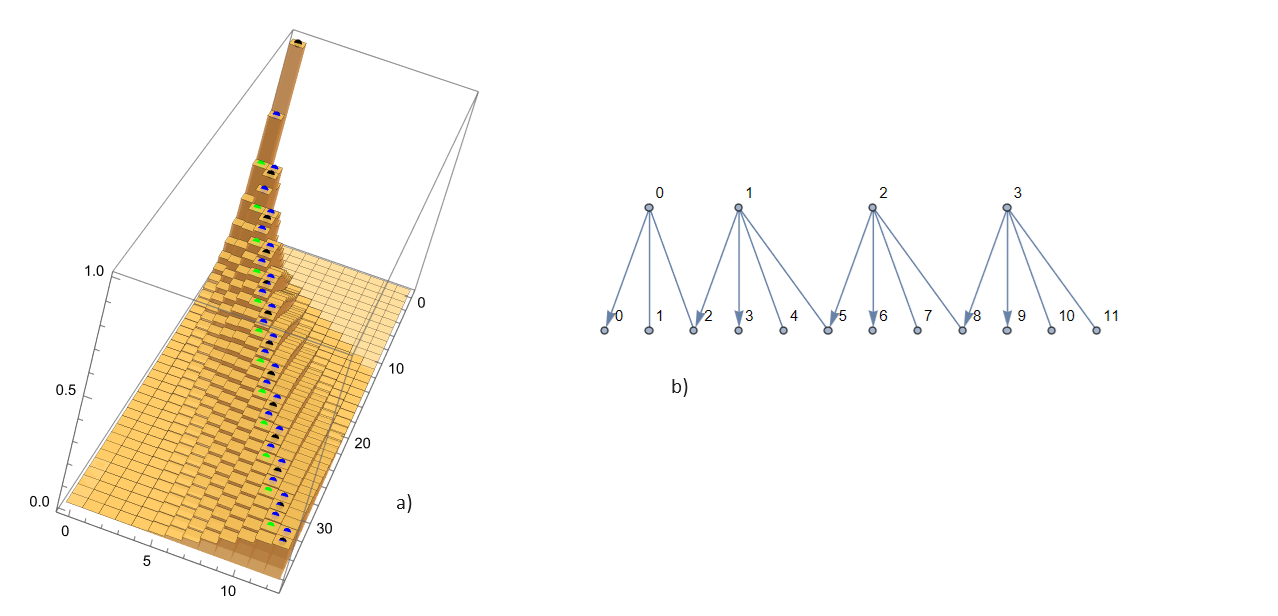}
\vspace*{-1.0cm}
\caption{a) The modal ridge of the binomial distribution with $p=1/3$ and varying $n$.
Blue and black spots mark the leading modes, black spots mark the larger likelihood maximiser.
b) 
For this distribution, the bipartite graph of the modal correspondence: e.g. $m_{\pm}(8,1/3)=\{2,3\}, \ell_\pm(2)=\{5,6\}$. 
 }
\label{RIDGE}
\end{figure}

For discrete distributions the situation is radically different, since the peak of probability function need not be constant exactly, rather has some room to vary along the modal ridge,
for instance to decay with $n$
as is the case for $S_n$ with binomial distribution, see Figure \ref{RIDGE}.

In this paper we focus on the modal properties of distributions representable as convolution of a Poisson distribution and any number (finite or  infinite)  of  Bernoulli distributions
with  success probabilities 
 that vary from trial to trial. 
Random variables with such distribution will be called here  {\it extended Bernoulli sums}.
This name echoes the term {\it extended Poisson-binomial distribution} suggested in 
 Broderick et al \cite{Broderick}.
In Madiman et al \cite{Madiman} the class of extended Bernoulli sums is introduced implicitly as a weak closure of the family of Poisson-binomial distributions.

The interest to  independent Bernoulli trials with varying success probabilities $p_1, p_2,\ldots$  is rooted deeply in history.
This sampling model was introduced by Poisson and 
 was known as the Poisson scheme in earlier literature,
but later the terminology was abandoned. 
The derived distribution of the number   of successes $S_n$ 
is 
nowadays commonly called  the  {\it Poisson-binomial} distribution \cite{Tang}. 
The first systematic study was undertaken in 1846 by
Chebyshev \cite{Chebyshev}, who established that amongst the Poisson-binomial distributions with given mean
the expectation ${\mathbb E}[\varphi(S_n)]$ is maximised by a shifted binomial distribution
(his main result was only stated for the
indicators needed to estimate the tail probabilities ${\mathbb P}[S_n\leq a]$ and  ${\mathbb P}[S_n\geq b]$
but the argument covers the general case).
About a century later  Hoeffding \cite{Hoeffding} proved a more detailed theorem
implying  that  the extremal property of the binomial distribution holds 
for  $\varphi$  any convex function on integers.
Bonferroni's  little known  paper \cite{Bonferroni} was the first to  observe the logconcavity of $S_n$, to  analyse  the 
behaviour of  the peak as $n$ varies and to show that the difference between the mode and the mean is always $o(n)$.
Apparently independentlly, the logconcavity  was rediscovered  soon thereafter by L{\'e}vy \cite{Levy}.
Darroch \cite{Darroch} and   Samuels \cite{Samuels}
found that  the mode differs from the mean by no more than $1$,  with the two being equal
 whenever the mean is an integer. Similar results for the median were shown in Jogdeo and Samuels \cite{Jogdeo}.
Pollard \cite{Pollard} presents a detailed exposition of much of this development. 
Pitman \cite{Pitman}  connected  the results about the mode  to  bounds on coefficients of  totally positive polynomials, giving many examples of combinatorial origin.
The modal properties of $S_n$ were most thoroughly studied  for the Karamata-Stirling profile of success probabilities $p_i=\theta/(\theta+i-1)$,
especially in the case $\theta=1$,
 in connection with Stirling numbers, theory of records and random combinatorial structures, 
see Kabluchko et al \cite{Kabluchko} for recent results and references.

In the other direction, the  statistical problem of estimating an unknown number of trials from the observed number of succeses  calls to consider $n$ as a parameter.
Assuming  the success probabilities known, 
the unimodality of likelihood  
 was  first  shown for  $p_i=1/i$  
 by Moreno-Rebollo et al \cite{Moreno1}  and  Cramer \cite{Cramer} then in general by Moreno-Rebollo et al \cite{Moreno2}. 
In disguise,  the fact also appeared in Bruss and Paindaveine \cite{Bruss} and  Tamaki \cite{Tamaki}  in connection with  versions of the famous secretary problem of optimal stopping.

The plan for  the rest of the paper is the following. In Section \ref{S2} we give a thorough review of  basic modal features of the Poisson and binomial distributions.
In Section \ref{S3} we show that for the logconcave power series families the cross modality amounts to Darroch's mean-mode rule.
In Section \ref{S4} we introduce techniques to work with extended Bernoulli sums and apply these to the analysis of the peak height.
Section \ref{S5} is devoted to the geometry of likelihood contours in the parameter space and presents an extension of Samuels \cite{Samuels} version of Darroch's rule.
In section \ref{S6} we prove the cross modality for certain monotonic (directed) families of extended Bernoulli sums, and give yet another proof of extended Darroch's rule by deriving it from the the cross modality.  
In the last section \ref{S7} we address aspects of a mass transport problem for the extended Bernoulli sums, 
with the objective to shift the mode.

In the sequel we will be dealing with families of discrete distributions, with context-dependent range $R$ of a family being either $\{0,1,\ldots,n\}$ or $\{0,1,\ldots\}$.

\section{The binomial and Poisson distributions}\label{S2}
In this  section  we revise under the angle of cross modality the well known properties of the binomial and the Poisson  distributions,
also using this occasion to introduce notation and terminology.
The connection between modes of the binomial and negative binomial distributions might have not been noticed in the literature.


For independent Bernoulli trials $B_1, B_2,\ldots$ with the same success probability $0<p<1$, let $S_n=B_1+\cdots+B_n$, $S_0=0$. 
The binomial distribution of $S_n$  has the
 probability function
\begin{equation}\label{bin-dis}
f(k; n,p) = {n\choose k} p^k (1-p)^{n-k},
\end{equation}
which we extend by zero outside the support $\{0,1,\ldots,n\}$.
The mode of this distribution is either a single integer or two adjacent integers (twins). 
The bifurcation of the mode can be neatly described by means of two  functions  
$m_{-}(n,p)\leq m_{+}(n,p)$,  where   $m_{+}(n,p)$ is   the {\it leading mode} determined  as the unique integer $k$
satisfying the inequalities
\begin{equation}\label{ineq-mode1}
f(k-1; n,p)\leq f(k; n,p)> f(k+1; n,p).
\end{equation}
The second function  $m_{-}(n,p)$ is defined by setting 
  $m_{-}(n,p)= m_{+}(n,p)$ unless
 the first relation in (\ref{ineq-mode1}) is equality in which case we set $m_{-}(n,p)= m_{+}(n,p)-1$.
With this notation, the mode of the binomial distribution is $m(n,p)=[m_-(n,p), m_+(n,p)]$, 
either singleton or two-point integer interval.
From the general perspective, the way of finding the mode via (\ref{ineq-mode1}) relies on the logconcavity of distribution.

Explicitly from (\ref{bin-dis})  and (\ref{ineq-mode1}),
 $m_+(n,p)$ is the unique integer $k$ satisfying  the inequalities
\begin{equation}\label{k-from-n}
k\leq (n+1)p< k+1. 
\end{equation}
If $(n+1)p$ is not integer then $m_{\pm}(n,p)=\lfloor(n+1)p\rfloor$ is a sole mode of the distribution, and if $(n+1)p$ is a positive integer then the mode 
is twin  $m(n,p)=[(n+1)p-1,(n+1)p]$.

To achieve solvability of the unconstrained likelihood maximisation problem we  close  the family of binomial distributions by allowing $p\in\{0,1\}$, and
extend the definition of modal functions as
$m_\pm(n,0)=0, m_\pm(n,1)=n$.
With the range of the family being $R=\{0,1,\ldots\}$,
the modal ridge can be identified  as  the set of triples $(k,n,p)$ where $k\in m(n,p), n\in\{0,1,\ldots\}$ and $p\in [0,1]$. 
The closed family is trivially cross modal, since for given $k$ the unique unconstrained likelihood maximiser is $(n, p)=(k,1)$ with $f(k;k,1)=1$.

We procede with  discussing the cross modality for two standard subfamilies  of binomial distributions, where the likelihood maximisers are genuine distributions.

\subsection{The binomial-$n$ family}
Suppose first that  { $0<p<1$ is fixed while $n$ varies}, 
thus  the likelihood becomes  a function  of the active parameter $n$. 
 Let 
$T_k$ be the index of  the $k$th success (with $T_0=0$). By the elementary identity
$$
{\mathbb P}[T_{k+1}=n+1]= p\, {\mathbb P}[S_{n}=k],
$$
the likelihood function $f(k; \cdot,p)$ is just a multiple of the probability function of the random variable $T_{k+1}-1$, which has a Pascal distribution. The distribution of  $T_k$ 
is the $k$th convolution power  of the geometric distribution for $T_1$ hence logconcave with one or two modes.
Similarly  to (\ref{ineq-mode1}) maximising the likelihood amounts to finding $\ell_+(n)$ as the unique $k$ satisfying
\begin{equation}\label{ineq-mode2}
f(k; n-1,p)\leq f(k; n,p)> f(k; n+1,p),
\end{equation}
then setting $\ell_-(n)=\ell_+(n)$ if the first relation in (\ref{ineq-mode2}) is strict and $\ell_-(n)=\ell_+(n)-1$ in the case of equality.
This leads to identifying $\ell_+(k)$ as the unique integer $n$ satisfying the inequalities
\begin{equation}\label{n-from-k}
n\leq \frac{k}{p}< n+1
\end{equation}
similar to (\ref{k-from-n}). Explicitly, if $\frac{k}{p}$ is not a positive integer we have   $\ell_+(n)=\ell_-(n)=\lfloor\frac{k}{p}  \rfloor$, and if $\frac{k}{p}$ is a positive integer
then $\ell_+(n)=\frac{k}{p}$ and $\ell_-(n)=\frac{k}{p}  -1$.

With the above explicit formulas in hand, it is now easy to verify the crossmodality of  the binomial sum $S_n$ relative to parameter $n$,
in the form
\begin{equation}\label{cr-m-bin}
n\in\{\ell_\pm(k)\}\Rightarrow k\in \{m_\pm(n)\},
\end{equation}
meaning that if $n$ is a likelihood maximiser for some $k$, then $(k,n)$ is a cross mode.

Indeed, suppose first that $\frac{k}{p}$ is not integer, thus  the inequalities for $n$ in
 (\ref{n-from-k}) are strict: then $n=\ell_\pm (k)$ and noting that  (\ref{n-from-k}) gives
$k<(n+1)p=np+p<np+1<k+1$
we obtain  (\ref{k-from-n}), whence
$k=m_{\pm}(n,p)$.
In the case $k=0$ we  trivially have   $m_{\pm}(0,p)=\ell_\pm(n)=0.$ 
In the remaining case of integer $n=kp$, we have that $n=\frac{k}{p}<n+1$ implies $\ell_+(k)=n$ and $\ell_-(k)=n=1$, that $k=np<k+1$ implies $m_+(n-1,p)=k$,  and $m_-(n-1,p)=k=1$,
and finally that $k<(n+1)p<k+1$ implies  $m_\pm (n,p)=k$ as above.
We conclude that  $m_{-}(\ell_-(k),p)=k-1$ only holds if  $\frac{k}{p}$ is a positive integer, while  $m_{\pm}(\ell_\pm(k),p)=k$ holds in all other cases.
See Figure \ref{RIDGE}.

Equivalently,  if $np<k<(n+1)p$ the pattern of relations about the cross mode $(k,n)$ is
\begin{equation}\label{c-m-ineq1}
\raisebox{.9\baselineskip}
{$_{\begin{array}{ccccc}
 &  & f(k;n-1,p) &  & \\
{} & {} & \vmore & {} & {} \\
 f (k-1;n,p)  & < &  f(k;n,p)  & >&  f(k+1;n,p)  \\
 & {} & \vless & {} & {} \\
 & {} &  f(k;n+1,p)  & {} & {} \end{array} }$}
\end{equation}
and if $k=(n+1)p$ the mode bifurcates and the pattern becomes 

\begin{equation}\label{c-m-ineq2}
\raisebox{.9\baselineskip}
{$_{\begin{array}{ccccc}
 &  & f(k;n-1,p) &  & \\
{} & {} & \vmore & {} & {} \\
 f (k-1;n,p)  & = &  f(k;n,p)  & >&  f(k+1;n,p)  \\
 & {} & \veq & {} & {} \\
 & {} &  f(k;n+1,p)  & {} & {} \end{array} }$}
\end{equation}
If $p$ is irrational number the mode does not bifurcate, hence only the strict pattern (\ref{c-m-ineq1}) occurs. If $p$ is rational, we may write $p=\frac{a}{b}$ with mutually prime integers $a<b$, 
and  the triple equality in (\ref{c-m-ineq1}) simplifies as the elementary identity for binomial coefficients,
$$
{cb -1\choose ca-1}b(b-a)={cb -1\choose ca}ab = {cb \choose ca} a(b-a).
$$

 Explicit formulas we used  obscur the mechanics of cross modality,  which becomes  transparent by turning to the recursion 
\begin{equation}\label{rec-bin-main}
f(k; n+1,p)=p f(k-1; n,p)+(1-p) f(k; n,p),
\end{equation}
that obtains by either  conditioning $S_{n+1}$ on the outcome of one trial or, what is  equivalent, by noting that adding to $S_n$  an independent Bernoulli variable yields $S_{n+1}$.
Geometrically, the recursion means that the probability function of $S_{n+1}$ is a convex combination of the probability functions of $S_n$ and $S_n+1$, see Figure.
If $m_-(n,p)=m_+(n,p)$,  the shapes of distributions of $S_n$ and $S_{n+1}$ have discordant directions
of strict monotonicity only 
within $[m_+(n,p), m_+(n,p)+1]$, and if  $m_-(n,p)<m_+(n,p)$ then disagreement occurs  within $[m_+(n,p)-1, m_+(n,p)+1]$ where one of the shapes is flat.
It follows that the unimodality of $S_n$ implies that $S_{n+1}$ is unimodal with a lower peak, unless the mode of $S_n$ is twin in which case the peak remains the same
and  the new mode becomes $m_\pm(n+1,p)=m_+(n,p)>m_-(n,p)$. This way we connect the maximisers as
\begin{equation}\label{mpr1}
\ell_-(k)=\min\{n: m_{+}(n,p)=k\},
\end{equation}
and observe  that $f(k; n,p)$ strictly increases in $n$ as long as  the leading mode does not hit $k$, then after possibly staying flat one step strictly decreases.

For the modal ridge we have therefore  dual representations  via the $n$- and $k$-sections as
$$M_{(\cdot,p)}=\bigcup\limits_{n=0}^\infty  m(n,p)\times{\{n\}}=\bigcup_{k=0}^\infty \{k\}\times [\ell_-(k), \ell_+(k+1)-1].$$

\subsection{The binomial-$p$ family}
The second, more common setting, is where
 $n>0$ is fixed and $p$ varies.  Let $\ell_{\circ}$  denote the likelihood maximiser relative to the `unknown $p$'  family of the binomial distributions.
Equating the derivative
$$
\partial_p f(k;n,p)=\frac{k-np}{p(1-p)}f(k;n,p)
$$
to zero we find   that the maximiser  is single-valued and given by the textbook formula
$\ell_\circ(k)=\frac{k}{n}$.
The cross modality follows from the identity
$$m_\pm\left(n,\frac{k}{n}\right)=k,~~~0\leq k\leq n,$$
which is trivial for $k\in\{0,n\}$ and    for  $0<k<n$  
is a  specialisation of (\ref{k-from-n})  as
$$k<(n+1)\frac{k}{n}<k+1.$$

The function $m_+(n,p)$ is right-continuous, having
unit jumps at points  $p_k=\frac{k}{n+1}$ that satisfy
$$f(k-1;n,p)=f(k;n,p), ~0\leq k\leq n,$$
and can be also characterised as maximisers of $f(k+1; n,p)$ as is seen from the formula
\begin{equation}\label{rec-bin-main-diff}
\partial_p f(k; n+1,p)= f(k-1; n,p)- f(k; n,p),
\end{equation}
derived from (\ref{rec-bin-main}).
For  $p_k$ with $0<k<n$  the mode is twin. Note the interlacing pattern of  pivot points 
$$0=p_0=\ell_\circ(0)<p_1<\ell_\circ(1)<p_2<\cdots<\ell_\circ(n-1)<p_n<\ell_\circ(n)=1.$$

The cross modality comes forward as the following pattern: 
the likelihood function  $f(k;p,n)$ increases in $p$ on $[0,p_k]$ until $k$ becomes a mode,  
then keeps growing until reaching  the maximum at $\ell_\circ(k)$, then decays with 
 the modal status  lost at  $p_{k+1}$.
The modal ridge becomes (with $p_{n+1}:=1$)
$$
M_{(n,\cdot)}= \bigcup_{k=0}^n \{k\}\times [p_k, p_{k+1}],
$$

\subsection{The Poisson distribution}

\begin{figure}
\centering
\includegraphics[scale=0.4]{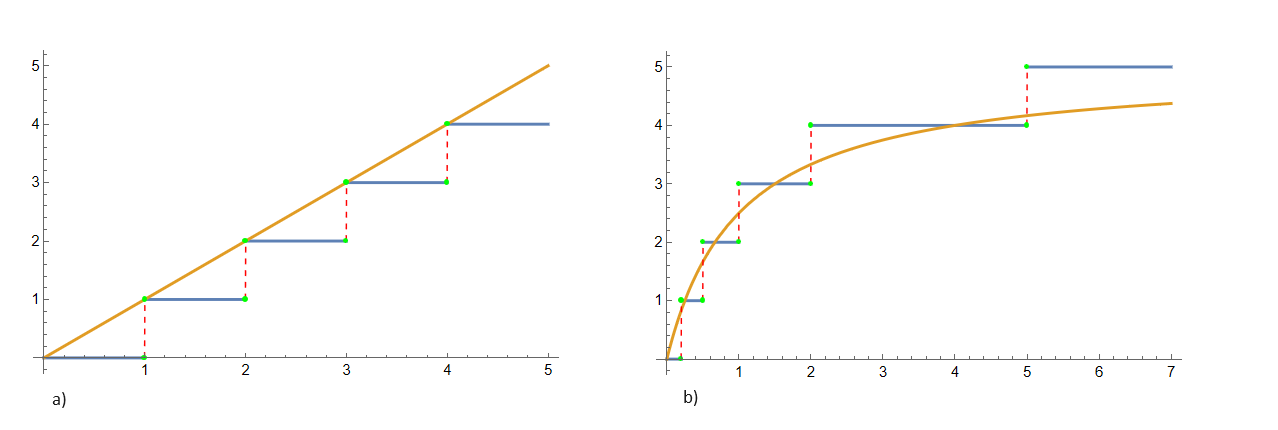}
\caption{a) The mean and the mode curves for the Poisson distribution. b) Same for the binomial distribution with $n=5$ and parameter $t=p/(1-p)$.
 }
\label{BinPo}
\end{figure}

For the Poisson distribution
$$
f(k; t)=e^{-t} \frac{t^k}{k!}, ~~~k\geq0,
$$
with parameter  $t\geq 0$ the modal properties are similar to that of both binomial families.
The pivotal values of the parameter are nonnegative integers $t_k=k$, whose role is twofold:

\begin{itemize}
\item[(i)] $f(k; t_k)=f(k-1; t_k), k\geq 1,$ hence the mode bifurcates, that is $m_+(t_k)=m_-(t_k)+1$, while $m_+(t)=m_-(t)=\lfloor t\rfloor$ for $t$ not a positive integer,

\item[(ii)] $f(k;t)$ as the likelihood  function coincides with the density of the Erlang distribution with shape parameter $k+1$, attaining its unique maximum at $t_k$. 
\end{itemize}
The likelihood maximisers coincide with the bifurcation points, forcing the mean curve to cross the graph of the leading mode in the northwest corner points, 
see the plot in Figure \ref{BinPo}.

The peak height $f(k; k)= e^{-k} \frac{k^k}{k!}$ is decreasing to $0$ with $k$.

\section{Power series distributions}\label{S3}

Every distribution on nonnegative integers can be included in one-parameter discrete exponential family, which can be conveniently represented by a power series with nonnegative coefficients.
A nice feature of such a family is that the likelihood maximisation amounts to  matching the mean with the observed value $k$.
In the unimodal case this allows one to express the cross modality in a simple way described in this section.

Let   $a_0, a_1,\ldots$ be a positive sequence such that the quotients of adjacent terms satisfy
$$\frac{a_{k+1}}{a_k} \downarrow 0~~{\rm as}~k\to\infty,$$
where the decrease is strict.
By the ratio convergence test,
the power series 
$$F(t)=\sum_{k=0}^\infty  a_k t^k$$
defines an entire function.
The associated power series distribution has the form
\begin{equation}\label{psd}
f(k; t)= \frac{a_k t^k}{F(t)},~~~k\geq0,
\end{equation}
where $t\geq0$ appears as a parameter. 
The probability generating function of (\ref{psd}) is 
$$z\mapsto \sum_{k=0}^\infty f(k;t)z^k=\frac{F(tz)}{F(t)}.$$

The assumption on coefficients implies strict logconcavity $a_k^2> a_{k-1} a_{k+1}$ for $k\geq 1$, hence the distribution is unimodal
and has all moments \cite{Keilson}.
The general  power series distribution 
 strictly increases with $t$ in the stochastic order, because the tail probability
 $\sum_{j=1}^\infty f(k+j;t)$ increases in consequence of  the fact that in the formula for the reciprocal
$$
\frac{1}{\sum_{j=1}^\infty f(k+j;t)}=
\frac{F(t)}{\sum_{j=1}^\infty a_{k+j} t^{k+j}}   = 1+   \sum_{i=0}^k \frac{a_i} {\sum_{j=1}^\infty a_{k+j} t^{k+j-i}}    $$
the terms of the external sum are decreasing to $0$. Therefore the limit distribution is 
the delta mass at $\infty$ as $t\to\infty$, and 
the mean
$$
\mu(t)= t\, \partial_t \log F(t) 
$$
is a smooth, strictly increasing function  with $\mu(0)=0$ and $\mu(t)\uparrow\infty$.

The likelihood quotient
$$\frac{f(k;t)}{f(k-1;t)}= \frac{a_{k}\,t}{a_{k-1}}, ~~~k\geq1,$$
is increasing  in both $k$ and $t$, and equals  $1$ at point
$$t_k:= \frac{a_{k-1}}{a_k},~~~k\geq1, $$
which is a unique solution to  $f(k-1;t)=f(k,t)$.
It is readily seen that $t_k\uparrow\infty$ (strictly, with $t_0:=0$), and that the leading mode $m_+(t)$ is right-continuous, nondecreasing,
takes value $k$ on the interval $[t_{k-1}, t_k]$ and bifurcates at the point $t_k$, where it has a unit jump for $k>0$ (that is, $m_-(t_k)=k-1$)).

We assert that  for each $k$ the likelihood has a single maximum point $\ell(k)$, which coincides with the (unique) solution to the equation $\mu(t)=k$.
Indeed, the likelihood is zero at $t\in\{0,\infty\}$,  and  the loglikelihood satisfies
$t\,\partial_t \log f(t;k)=k -\mu(t)$ by easy calculus.

\begin{prop} \label{prop1}
The    family of logconcave power series distributions
{\rm (\ref{psd})} 
is cross modal if and only if whenever $\mu(t)$ is integer it coincides with a mode.
This condition is equivalent to any of the following three conditions:
\begin{itemize}
\item[{\rm (i)}] $t_{k}\leq \ell(k)\leq t_{k+1},~~k\geq 1$,
\item[{\rm (ii)}] $ k-1\leq \mu(t_k)\leq k, ~~k\geq 1$,
\item[{\rm (iii)}] $\sup |\mu(t)-m_+(t)|\leq 1$.
\end{itemize}
\end{prop}

\proof  For cross modality  $k$ must be a mode of $f(\cdot;\ell (k))$, which is (i) since $k$ is a mode for $t\in[t_{k-1}, t_k]$.
The rest of proposition follows from the identification of $\ell(k)$ as the root of $\mu(t)=k$, combined with the  continuity and monotonicity properties of $m_+(t)$ and $\mu(t)$.
\endpf

\begin{figure}
\centering
\includegraphics[scale=0.5]{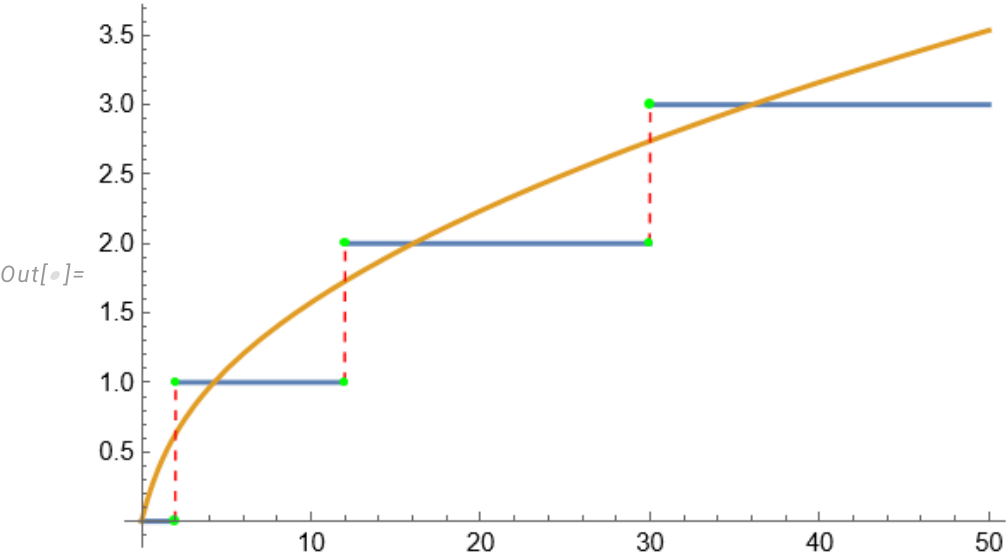}
\caption{The mean and the mode functions for the power series distribution derived from $F(t)=\cosh \sqrt{t}$. }
\label{COSH}
\end{figure}

Explicitly, condition (ii) becomes
\begin{equation}\label{cmd-odd}
(k-1)\sum_{j=0}^\infty a_j\left( \frac{a_{k-1}}{a_k}\right)^j \leq \sum_{j=1}^\infty j a_j\left( \frac{a_{k-1}}{a_k}\right)^j\leq k\sum_{j=0}^\infty a_j\left( \frac{a_{k-1}}{a_k}\right)^j.
\end{equation}
Since the Poisson distribution is the edge case, with $\ell(t_k)=\mu(t_k)=t_k=k$ and the second relation in (\ref{cmd-odd}) being equality, one can hope that the ultra logconcavity
$$
\frac{a_k}{a_{k+1}}\frac{1}{k+1}\geq \frac{a_{k-1}}{a_{k}}\frac{1}{k}
$$
would be sufficient for the cross modality. Unfortunately, this does not materialise. For,  denote the right-hand side $b_k:= a_{k-1}/(k a_k)$; these parameters are now only constrained by the monotonicity $0<b_1\leq b_2\leq\ldots$
Manipulating the second inequality in (\ref{cmd-odd}) we reduce it to the condition
$$ 
\sum_{j=0}^\infty \frac{b_k^j}{b_1\cdots b_j} \frac{k^j}{j!}\left( 1- \frac{b_k}{b_{j+1}} \right) \geq 0.
$$
But for $k>1$ and any given $b_k, b_{k+1},\ldots$ we may choose $b_1<\ldots<b_{k-1}$ sufficiently small to have the whole infinite sum negative.

A finite positive sequence $a_0,\ldots,a_n$ seen as coefficients of the polynomial $F(t)=\sum_{k=1}^n a_k t^k$ induces a family of distributions supported by $\{0,1,\ldots,n\}$.
The definitions and findings of this section are easily adjusted to the polynomial case.
In this framework the largest bifurcation point 
is $t_n=a_{n-1}/a_n$ for polynomial of degree $n$, and the mean satisfies $\mu(t)\uparrow n$, so it is logical to set $\ell(n)=\infty$.
The  binomial-$p$ family fits in the scheme by the re-parametrisation $t=p/(1-p)$, see Figure \ref{BinPo}.
Pitman \cite{Pitman} argues in the direction opposite to ours, observing that the mean-mode
relations imply bounds on the quotients $a_{k-1}/a_{k}$.

\section{Extended binomial sums}\label{S4}

\subsection{Definitions and construction}

Given $\lambda\geq 0$ and a profile of success probabilities ${\bp}=(p_1, p_2,\ldots)$ with  $\sum_{i=1}^\infty  p_i  <\infty$,
 the extended Bernoulli sum is a random variable
$S=X+\sum_{i=1}^\infty B_i,$ where $X\stackrel{d}{=}{\rm Poisson}(\lambda)$ and $B_i\stackrel{d}{=}{\rm Bernoulli}(p_i)$ are jointly  independent.
The probability function of $S$, written as $f(k; \lambda,\bp),~k\geq0$, 
depends symmetrically on the $p_i$'s, so  the uniqueness of parametrisation is  achieved by
relating each distribution to a point of the infinite-dimensional convex set
$\nabla:=\{(\lambda,\bp):\lambda\geq 0~{\rm and~} 1\geq p_1\geq p_2\geq \ldots\geq 0)$, which 
has  the natural  component-wise partial order. 
The mean of $S$ will be denoted 
 $$
\mu(\lambda,\bp):=\sum_{k=0}^\infty k f(k;\lambda,\bp)=\lambda+\sum_{i=1}^\infty p_i.
$$
Note that $\mu(\lambda,\bp)$ is the $l_1$-norm, and
whenever  $(\lambda_1,\bp_1)\geq (\lambda_2,\bp_2)$ the difference  $\mu(\lambda_1,\bp_1)-\mu(\lambda_2,\bp_2)$ equals the $l_1$-distance between the vectors $(\lambda_1,\bp_1)$ and  $(\lambda_2,\bp_2)$.

We endow $\nabla$ with the topology in which convergence  amounts to the  convergence in distribution of the corresponding extended Bernoulli sums
(that is, pointwise convergence of the $f(\cdot;\lambda,\bp)$'s).
Explicitly,  convergence of  a sequence $\{(\lambda_j, \bp_j)\}\subset\nabla$
in this topology amounts to the componentwise convergence of $\bp_j$'s {\it and} convergence of the mean values $\mu(\lambda_j,\bp_j)$,
with the behaviour of $\lambda_j$'s being subordinate to these rules.
In particular,   the familiar convergence of Bernoulli convolutions to a Poisson distribution in terms of the parameters is equivalent to
$(0, \bp_j)\to (\lambda,{\boldsymbol 0})$,
which holds when
$\bp_j\to {\boldsymbol 0}$ (this is equivalent to the convergence of the largest term to $0$)
and  $\mu(0,\bp_j)\to \lambda$.
To illustrate,  convergence of the  binomial distribution to Poisson corresponds to the relation
$$(0,\underbrace{\lambda/n,\ldots,\lambda/n}_ {n~ {\rm times}}, 0,0,\ldots) \to (\lambda,0,0,\ldots).$$

By the alternative parametrisation by  $\mu$ and ordered  and ordered $\bp$ satisfying $\mu\geq \sum_{i=1}^\infty p_i$, the adopted topology  becomes the product topology.
Broderick et al \cite{Broderick} use  discrete measures $\lambda\delta_0+ \sum_{i=1}^\infty p_i \delta_{p_i}$
to encode the parameters, but the vector representation   has the advantage that the Bernoulli terms get naturally labelled.
The uniqueness of parametrisation by $\nabla$ is a part of a more general result on totally positive sequences \cite{Karlin}. See \cite{Broderick} for a more probabilistic proof.

The {\it finitary} Poisson-binomial distributions with   
$n$  Bernoulli terms (some of which could be $0$)
correspond to a  compact face  $\nabla_n\subset\nabla$ with points of the form $(0,p_1,\ldots,p_n,0,0,\ldots)$. The union $\cup_n \nabla_n$ is dense in $\nabla$,
therefore many structural properties of the extended Poisson-binomial distribution can be concluded from their 
finitary counterparts
but we prefer to avoid such way of reason, giving  direct proofs whenever possible.

A simple  randomisation device 
allows one to organise all extended Bernoulli sums in a infinite-parameter counting process  $(S(\lambda,\bp),    ~(\lambda,\bp)\in {\mathbb R}_+\times [0,1])$, where we do not restrict $\bp$ to have the  components 
$p_i$
nondecreasing.
 To that end, let $(N(\lambda), \lambda\geq 0)$ be a unit rate Poisson process,
 independent of   i.i.d.  random variables $U_1, U_2,\ldots$ with uniform distribution on $[0,1]$;
this data can be regarded as a point process on the ground space ${\mathbb R}_+\cup (\cup_{i=1}^\infty [0,1])$.
Setting $S(\lambda,\bp)=N(\lambda)+ \sum_{i=1}^\infty 1(U_i\leq p_i)$ we obtain an integer-valued process with desired 
marginal distributions  and paths nondecreasing in $(\lambda,\bp)$. Then $S(\lambda,\bp)$ is almost surely finite if $\mu(\lambda,\bp)<\infty$.

In the sequel, unless explicitly required, we will not restrict $p_i$'s to be nondecreasing. This freedom is needed to avoid re-labelling while allowing a particular success probability  parameter to vary within the full range $[0,1$.

Sometimes it is convenient to use a parallel set of  parameters $\br=(r_1, r_2,\ldots)$ for the Bernoulli variables with the largest $p_1<1$, where $r_i$ is
interpreted as odds in favour of success, and the connection is 
$$r_i:=\frac{p_i}{1-p_i}\Longleftrightarrow p_i=\frac{r_i}{1+r_i}.$$
Clearly, $\sum_i r_i<\infty$ is equivalent to $\sum_i p_i<\infty$ if $p_1<1$.

The probability generating function of $S$ 
is an entire function with only negative zeros, 
\begin{equation}\label{entireGF}
G(z)=G(0)\, e^{\lambda z} \prod_{i=1}^\infty (1+r_i z), ~~{\rm where}~~G(0)= e^{-\lambda} \prod_{i=1}^\infty (1-p_i)= e^{-\lambda} \prod_{i=1}^\infty \frac{1}{1+r_i}.
\end{equation}
For $k=0,1,\ldots$ 
consider the elementary symmetric functions in the odds
$$
E_k(0,\br):=\sum_{i_1<\cdots<i_k} r_{i_1}\cdots r_{i_k},
$$
and let us introduce their extended counterparts
$$
E_k(\lambda,\br):=\sum_{j=0}^k \frac{\lambda^{k-j}}{(k-j)!} E_j(0,\br).
$$
The convergence of the multiple series is ensured by the summability of $r_i$'s 
and can be concluded by expressing the elementary symmetric functions as polynomials in Newton's power sums $\sum_i r_i^n$.
In this notation, the probability function of $S$ becomes
$$
f(0;\lambda, \bp)= S(0),~~~f(k;\lambda, \bp)=f(0;\lambda, \bp) E_k(\lambda,\br).
$$
We extend this by zero, setting $ f(-k; \lambda,\bp)=E(-k; \lambda,\bp)=0$ for $k>0$.

Throughout $f(k;\lambda,\bp\setminus p_i)$ will denote the probability of $k$ successes with the $i$th Bernoulli trial ignored
(this is the probability function of $S-B_i$),    and  $E_k(\lambda,\br \setminus r_i)$ will denote the symmetric function with variable $r_i$ set to zero.
Reciprocally, the convolution with yet another independent Bernoulli$(p)$ variable is denoted $f(k;\lambda,\bp\cup p)$.
Similar notation will apply to the operations of removing or adding two or more independent Bernoulli variables.
The operator $\Delta$ will  act in the variable $k$ as forward difference
$\Delta f(k)=f(k+1)-f(k)$.

The following lemma summarises some of the properties of the probability function.

\begin{lemma}\label{L1} The probability function of the extended Bernoulli sum satisfies:

\begin{itemize}
\item[\rm(i)] $f(k; \lambda,\bp)=p_i    f(k-1; \lambda,\bp\setminus p_i)+(1-p_i)    f(k; \lambda,\bp\setminus p_i),$

\item[\rm(ii)] {\rm (}ultra log-concavity{\rm )} for $k\geq1$

$$ {k f^2(k; \lambda,\bp)} \geq {(k+1){f(k-1; \lambda,\bp)}f(k+1; \lambda,\bp)} ,$$
where the inequality is strict unless $\bp=0$, or $\lambda=0$ and there are less than 
$k$ positive $p_i$'s {\rm(}when both sides are zero{\rm )},

\item[\rm(iii)] {\rm (}linearity in $p_i${\rm )}
$$\partial_{p_i} f(k; \lambda,\bp)=  -\Delta  f(k-1; \lambda,\bp\setminus p_i),$$
\item[\rm(iv)] $\partial_\lambda f(k; \lambda,\bp)=-\Delta     f(k-1; \lambda,\bp),$

\item[\rm(v)] {\rm (}monotone likelihood ratio{\rm )}
$$\partial_{p_i} \frac{f(k+1; \lambda,\bp)}{f(k; \lambda,\bp)}>0, ~~~\partial_{\lambda} \frac{f(k+1; \lambda,\bp)}{f(k; \lambda,\bp)}>0,$$

\item[\rm(vi)]    $k f(k;\lambda,\bp)=\sum_{i=1}^\infty p_i f(k;\lambda,\bp\setminus p_i)+\lambda f(k-1;\lambda,\bp).$

\item[\rm(vii)]   {\rm (}concave likelihood ratio{\rm )}
$$\frac{f(k+1;\lambda,\bp)}{f(k;\lambda,\bp)}$$
 is concave in $(\lambda,\bp)$, and strictly concave for $k>0$.

\item[\rm(viii)] for $d\geq 1$ and  pairwise distinct $i_1,\ldots,i_d$
$$\partial^d_{p_{i_1},\ldots,p_{i_d}} f(k;\lambda,\bp)=(-\Delta)^d f(k-d;\lambda,\bp\setminus p_{i_1},\ldots,p_{i_d}).$$

\end{itemize}
\end{lemma}
\proof  Assertion  (i) follows by the total probability rule, which in terms of the symmetric functions becomes the recursion
$$ E_{k}(\lambda,\br)=r_i E_{k-1}(\lambda,\br\setminus r_i)+E_{k}(\lambda,\br\setminus r_i).$$
Formula  (iii) follows straight from (i).

Part (ii) for $\bp=0$ is an easy factorial identity. For  some number $n\geq k$ of nonzero $p_i$'s, matching powers of $\lambda$ and applying Newton's inequality 
yields a stronger inequality. The case with infinitely many positive $p_i$'s  is treated by induction in $k$ using (i).
The claim 
also follows by checking the property for Poisson and Bernoulli distributions, and applying
 Liggett's theorem \cite{Liggett} which asserts that the ultra log-concavity is preserved by convolution,
Saumard and 
Wellner \cite{Wellner} 
give a detailed discussion, including
the  interpretation  of ultra-logconcavity 
as the log-concavity  relative to the Poisson distribution.

Part (iv) is straight from $\partial_\lambda f(0;\lambda,\bp)=- f(0;\lambda,\bp)$ and $\partial_\lambda E_k(\lambda,\br)=E_{k-1}(\lambda,\br)$.

For (v), using $\partial_{p_i}r_i>0$ we are reduced to checking the analogue in terms of the odds, which is the inequality
$\partial_{r_i}[E_{k+1}(\lambda,\br)/ E_{k}(\lambda,\br)]>0.$
Differentiating the quotient, after simplification of the numerator the inequality follows from log-concavity.

The identity (vi) was stated in \cite{Samuels} for the finitary case without proof.
To show (vi), write the desired identity in terms of the odds,
$$k E_k(\lambda,\br)=\sum_{i=1}^\infty r_i E_k(\lambda,\br\setminus r_i)+\lambda E_{k-1}(\lambda,\br).$$
The formula is proved by differentiating the homogeneity identity $E_k(z\lambda,z\br)=z^k E_k(\lambda,\br)$ in $z>0$, using the above  formulas for partial derivatives,  and setting
 $z=1$.

(vii) Passing to the quotient of $E_k(\lambda,\br)$'s this is shown as in \cite{Marshall}, p. 116.

(viii) Obtained by iterating (i) or (iii).

\endpf
\noindent
A number of other  classical inequalities for symmetric polynomials \cite{Marshall}
that do not restrict the number of variables
 have their extended analogues and can be interpreted 
as features of the probability function $f(k;\lambda,\bp)$.


By the log-concavity we can  define  the leading mode $m_+(\lambda,\bp)$ to be the unique integer $k$ satisfying 
 $$f(k-1,\lambda,\bp)\leq f(k,\lambda,\bp)>  f(k+1,\lambda,\bp),$$
and  define $m_-(\lambda,\bp)$ as in the case of the binomial distribution.
The probability function strictly increases for $k<m_-(\lambda,\bp)$ and strictly decreases for $k\geq m_+(\lambda,\bp)$ (unless turning zero in the finitary case). 
If $m_-(\lambda,\bp)=m_+(\lambda,\bp)$ the mode is singleton, and if  $m_-(\lambda,\bp)+1=m_+(\lambda,\bp)$ the mode is twin.
Thus, the mode is identified by the equivalences
\begin{eqnarray}
\label{E1}
f(k-1;\lambda,\bp)< f(k;\lambda,\bp) &\Longleftrightarrow&  m_- (\lambda,\bp)\geq k,\\
\label{E2}
f(k;\lambda,\bp)> f(k+1;\lambda,\bp)& \Longleftrightarrow & m_+ (\lambda,\bp)  \leq k.
\end{eqnarray}
The peak height of the distribution 
$h(\lambda,\bp):=\max_k f(k;\lambda,\bp)$
is attained at both $m_+(\lambda,\bp)$ and $m_-(\lambda,\bp)$ even  when they are distinct, but only the leading mode
 $m_+(\lambda,\bp)$ is stable  under  small {\it  increase} of the parameters.

\subsection{The peak height}

The peak height is a decreasing function of $\lambda$, that is $\partial_{\lambda}h(\lambda,\bp)\leq 0$ where $0$ only occurs if the   mode at $(\lambda,\bp)$ is twin.
We move on to revealing  a more involved  dependence of the peak on the $p_i$'s. 

\begin{figure}
\centering
\includegraphics[scale=0.5]{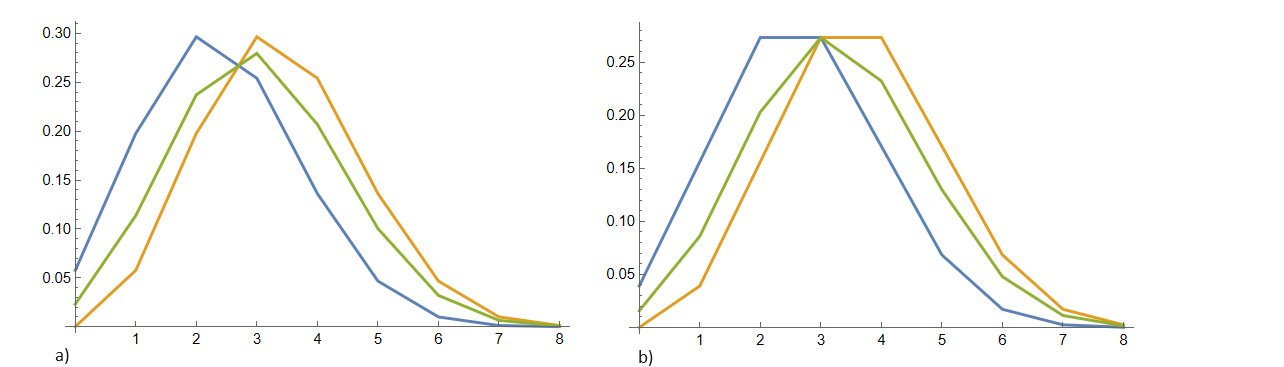}
\caption{a) Evolution of single mode by adding a Bernoulli variable. b) Same for twin. }
\vspace*{-0.5cm}
\label{SHIFT}
\end{figure}

To start with, recall  the discussion around (\ref{rec-bin-main}) for the binomial distribution.
Similarly to that, we may conclude on the
unimodality  of $S$  straight from the recursion  (i) in Lemma \ref{L1}, first using induction in the number of Bernoulli components for finitary distributions, then passing 
to limit.   The key step is the following.
Fix $(\lambda,\bp)$, and suppose yet another independent Bernoulli variable $B$ is added to $S$.
Thus we have as in (\ref{rec-bin-main})
\begin{equation}\label{C-C}
f(k;  \lambda, {\bp}\cup p)= p  \, f(k-1;  \lambda, {\bp}) +  (1-p)  f(k;  \lambda, {\bp}),
\end{equation}
where the right-hand side is to be understood as a convex combination of probability function and its unit shift.
The shapes in Figure \ref{SHIFT} should convince the reader that the unimodality of $f(k;  \lambda, {\bp}\cup p)$ follows from that of $f(k;  \lambda, {\bp})$.

A quick  look at Figure \ref{SHIFT} shows that the shapes of (linearly interpolated) probability functions in (\ref{C-C}) have a sole crossing point, which does not depend on $p$.
Write for shorthand $m_+$ for the leading mode, and   $f(m_+-1), f(m_+), f(m_++1)$ for
the values of  $f(\cdot;  \lambda, {\bp})$ at these positions.
The height of the crossing point is computed as 
$$ \frac{f^2(m_+) -f(m_+-1)f(m_++1)}{2f(m_+)- f(m_+-1)- fm_++1)}.$$
The peak height of $f(\cdot\,;  \lambda, {\bp}\cup p)$ will be at this level if $p$ is equal to the quantity 
$$
\gamma^*(\lambda,\bp):=\frac{-\Delta f(m_+) }{  \Delta^2 f(m_+-1)}\,,
$$
which we will call the {\it peak skewness} (not to be confused with distribution skewness). 
A similar characteristic of the peak called `ratio' was introduced in Baker and Handelman \cite{Baker}.
The leading mode of $S+B$ is $m_+$ if $p<\gamma^*(\lambda,\bp)$ and switches to $m_+$ if $p\geq \gamma^*(\lambda,\bp)$.
Peak skewness $1$ only occurs if the mode of $S$ is twin, while peak skewness $1/2$ means that the distribution is symmetric near the mode.
The peak height of $S+B$ obtains as 

$$
h(\lambda,\bp\cup p)=
\begin{cases}
(1-p)f(m_+)+p f(m_++1),~~~{\rm if~}p>\gamma^*(\lambda,\bp),\\
(1-p)f(m_+)+p f(m_+-1),~~~{\rm if~}p<\gamma^*(\lambda,\bp),\\
f(m_+), ~~~{\rm if~}p=\gamma^*(\lambda,\bp).
\end{cases}
$$
If the mode of $S$ is twin, the peak height of $S+B$ is the same as for $S$ regardless of $p$.

The above analysis of the effect from adding a Bernoulli variable readily tells us what happens when a variable gets removed from the extended Bernoulli sum.

\begin{theorem}\label{T1} Both $m_+(\lambda,\bp)$ and $m_-(\lambda,\bp)$ are nondecreasing functions of the parameters, with the following connection to the monotonicity of the peak.
The peak height function is differentiable in $\lambda$ and every $p_i$ and for $k= m_+ (\lambda,\bp)$  we have
\begin{itemize}
\item[\rm (i)] $\partial_{p_i}h(\lambda,\bp)>0$ iff $m_+(\lambda,\bp\setminus p_i)=m_-(\lambda,\bp\setminus p_i)=k-1,$ in which case replacing $p_i$  in $\bp$ by any larger value
$p_i+\delta$ cannot switch the leading mode to $k+1$, 
\item[\rm (ii)] $\partial_{p_i}h(\lambda,\bp)<0$ iff $m_+(\lambda,\bp\setminus p_i)=m_-(\lambda,\bp\setminus p_i)=k,$ in which case 
replacing $p_i$  in $\bp$ by certain larger value
$p_i+\delta$   will switch the leading mode to $k+1$,
\item[\rm (iii)] $\partial_{p_i}h(\lambda,\bp)=0$ iff $m_+(\lambda,\bp\setminus p_i)=m_-(\lambda,\bp\setminus p_i)+1=k,$ in which case replacing $p_i$  in $\bp$  
 by any larger value
$p_i+\delta$  will leave the single mode at $k$.
\item[\rm (iv)] $\partial_{\lambda}h(\lambda,\bp)\leq 0.$

\end{itemize}
\end{theorem}
\proof
We have seen that manipulating $p_i$ alone results 
in two possibilities for the mode, depending on how $p_i$ compares to the peak skewness of $f(\cdot;\lambda, \bp\setminus p_i)$. 
For instance, in case (i) replacing $p_i$ by $0$  moves the mode to  $k=1$, 
hence $k+1$ cannot be achieved by  replacing $p_i$ with $p_i+\delta$.
The sign of  $\partial_{p_i}h(\lambda,\bp)$ reads from Lemma \ref{L1}(iii) specialised for  $k$ the modal value.
Part (iv) follows from Lemma \ref{L1} (iv).
\endpf
\noindent
Whenever different monotonicity directions in (i) and (ii) occur the peak heaight function has a saddle point.

\section{Geometry of the  Darroch  rule}\label{S5}

\subsection{The Poisson-binomial case}

Darroch's \cite{Darroch} rule localises   the mode of the Poisson-binomial distribution in terms of the mean $\mu$.
The rule states that whenever the mean is integer it coincides with a single mode, and if the mean is not a positive integer then 
the mode could be singleton $\lfloor\mu\rfloor $ or $\lceil\mu\rceil $ , or twin $\{\lfloor\mu\rfloor, \lceil\mu\rceil \}.$
Samuels \cite{Samuels} has noticed that the probability function is increasing for $k\leq\mu$ and decreasing for $k\geq\mu$.

A loose intuition suggesting that the mode and the mean are not far away from one another appears by looking at two extremes.
If a Bernoulli variable with success probability zero or one is added to $S$ the mean and the mode both increment by zero or one, respectively.
It will take some effort to make this idea work.

We give a geometric interpretation to the rule first in the finitary case.
Consider  the $n$-dimensional simplex $\nabla_n\subset\nabla$ whose points are the vectors of parameters $(0,p_1,\ldots,p_n,0,0,\ldots)$, where 
$1\geq p_1\geq\ldots\geq p_n\geq 0$. (The first zero in this notation stays for $\lambda=0$.)
The apex of $\nabla_n$ is at $(0,{\boldsymbol 0})$ and the simplex itself is a cap of  a lattice cone.
Let 
\begin{eqnarray*}
{\cal F}_{k,n}&=&\{(0,\bp)\in\nabla_n: f(k-1; 0,\bp)=f(k; 0,\bp)\}, ~k\geq0,\\
{\cal M}_{x,n}&=&\{(0,\bp)\in\nabla_n: \mu(0,\bp)=x\},~x\geq0.
\end{eqnarray*}
For $0<k\leq n $,         ${\cal F}_{k,n}$ is a $(n-1)$-dimensional algebraic manifold (with boundary)
which can be regarded as a section
splitting $\nabla_n$ in two connected components. 
For $0<x<n$, ${\cal M}_{x,n}$ is a  section of  $\nabla_n$ by a hyperplane.
For integer $0<k<n$ the set ${\cal M}_{k,n}$ is a $(n-1)$ dimensional simplex, for instance ${\cal M}_{2,3}$ is a triangle with extreme points
$$(0,1,1,0,0,\ldots),   (0,1,\tfrac{1}{2},\tfrac{1}{2},0,\ldots), (0,\tfrac{2}{3}, \tfrac{2}{3},\tfrac{2}{3},0,0,\ldots).$$
The ${\cal F}_{k,n}$'s for $0\leq k\leq n$ are pairwise disjoint in accord with the fact that 
no three positive values of the   Poisson-binomial probability function can be equal.

\begin{figure}
\centering
\includegraphics[scale=0.5]{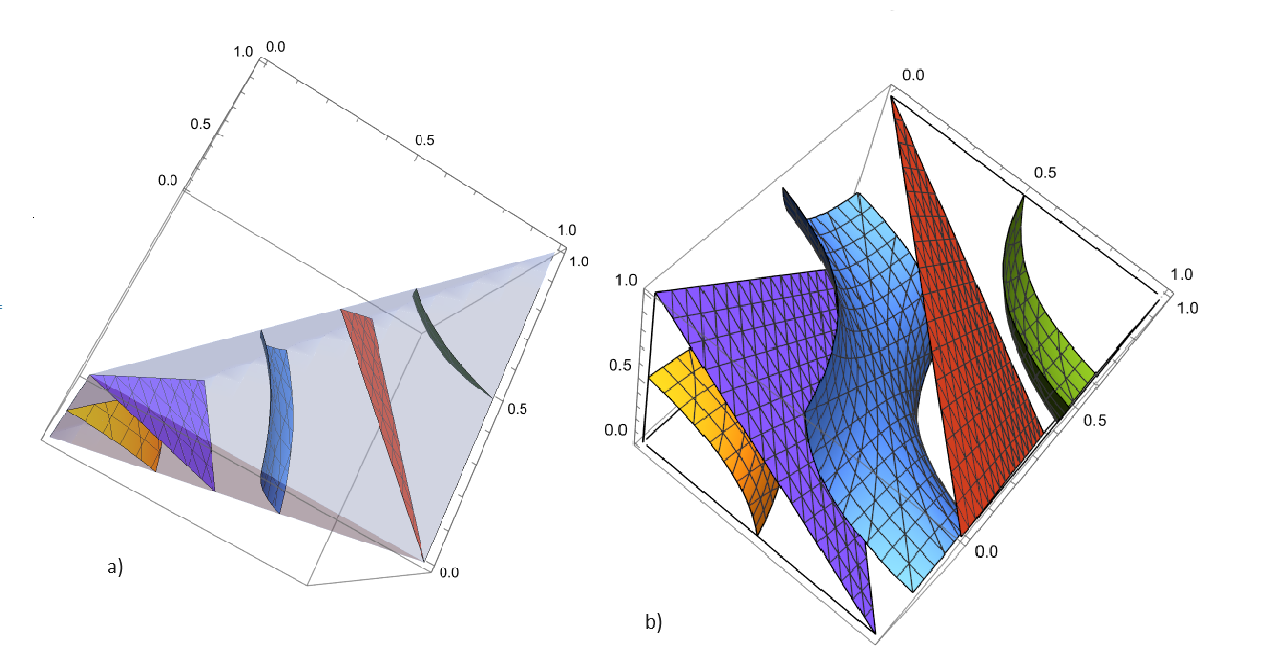}
\caption{For $n=3$ surfaces $f(0)=f(1), \mu=1, f(1)=f(2), \mu=2$ and  $f(2)=f(3)$ in a) simplex $\nabla_3$, and b) cube $[0,1]^3$.
Shorthand notation: $f(k)=f(k;0,p_1,p_2,p_3,0,0,\ldots)$, same for $\mu$. 
 }
\vspace*{-0.5cm}
\label{PBGeom1}
\end{figure}

By Darroch's rule the manifolds arranged in the interlacing sequence
$$  {\cal F}_{1,n}, {\cal M}_{1,n}, {\cal F}_{2,n},\ldots, {\cal M}_{n-1,n}, {\cal F}_{n,n}$$
are pairwise disjoint, as they are separated in terms of the $l_1$-norm
$\mu(0,\bp)$, see Figure \ref{PBGeom1}.
For two neighbours in the sequence we may identify  pairs of their points
most close and most distant in the $l_1$-metric. 
Explicitly, the unique pair
\begin{eqnarray*}
(0,\underbrace{1,\ldots,1}_{k-1}, \underbrace{\frac{1}{n-k+2},\ldots,\frac{1}{n-k+2}}_{n-k+1},0,0,\ldots)&\in& {\cal F}_{k,n},\\
(0,\underbrace{1,\ldots,1}_{k-1}, \underbrace{\frac{1}{n-k+1},\ldots,\frac{1}{n-k+1}}_{n-k+1},0,0,\ldots)&\in& {\cal M}_{k,n}
\end{eqnarray*}
realises the largest $l_1$-distance  between the sets, which is equal to  $ 1-\frac{n-k+1}{n-k+2}$.
Likewise, the unique pair 
\begin{eqnarray*}
(0,\underbrace{\frac{k}{k+1},\ldots,\frac{k}{k+1}}_{k},0,0,\ldots)&\in& {\cal F}_{k,n},\\
(0,\underbrace{1,\ldots,1}_{k},0,0,\ldots)&\in& {\cal M}_{k,n}
\end{eqnarray*}
realises the smallest $l_1$-distance  between the sets, which is  equal to $ 1-\frac{1}{k+1}$.
This can be verified using Chebyshev's \cite{Chebyshev} device that the extreme values of $f(k; 0,\bp)-f(k-1; 0,\bp)$ on ${\cal M}_{x,n}$ are attained at some
shifted binomial probability functions. An alternative approach is to reduce to a  nonlinear programming problem with concave functions $f(k)/f(k-1)$ (Lemma \ref{L1}  (vii)) and concave $\mu$ all expressed in the odds variables.
It follows that for $0<k\leq n$ and $x$ in the bounds 
\begin{equation}\label{maxmin}
k-1+\frac{1}{k+1}= \min_{{\cal F}_{k,n}}               \mu(0,\bp)    <x< \max_{{\cal F}_{k,n}}  \mu(0,\bp)=k-\frac{1}{n-k+2}\,,
\end{equation}
for $\bp\in {\cal M}_{x,n}$ the mode could be any of $k, k-1, \{k-1,k\}$ and each possibility materialises.
This leaves us for  $0\leq k\leq n$   with the definite region
\begin{eqnarray}
k-\frac{1}{n-k+2}<\mu(0,\bp)< k+\frac{1}{k+2}&\Longrightarrow & m(0,\bp)=\{k\}, \nonumber \\
\mu(0,\bp)=k-\frac{1}{n-k+2}&\Longrightarrow&  m(0,\bp)=\{k-1,k\}, \nonumber\\
\mu(0,\bp)= k+\frac{1}{k+2}&\Longrightarrow&  m(0,\bp)=\{k,k+1\},
\label{Dar}
\end{eqnarray}
which together with the assertion of ambiguity for the complementary region (\ref{maxmin}) constute a more detailed version of Darroch's rule.

\subsection{The extended  binomial sums}

We turn next to the extended binomial sums. Introduce
\begin{eqnarray*}
{\cal F}_{k}&=&\{(\lambda,\bp)\in\nabla: f(k-1; \lambda,\bp)=f(k; \lambda,\bp)\}, ~k\geq0,\\
{\cal M}_{x}&=&\{(\lambda,\bp)\in\nabla: \mu(\lambda,\bp)=x\},~x\geq0.
\end{eqnarray*}
Similarly to the finitary case, the ${\cal F}_{k}$'s, $k\geq0$, are pairwise disjoint closed sets
(`algebraic manifolds of codimension one'),
and the ${\cal M}_{x}$'s, $x\geq0$,  are pairwise disjoint convex  closed sets.
By the definition we have
$$ {\cal F}_{k,n}={\cal F}_{k}\cap \nabla_n, ~~~ {\cal M}_{k,n}={\cal M}_{k}\cap \nabla_n,$$
which suggests that an infinite counterpart of Darroch's rule could be obtained by an approximation argument.
This is true and not hard to pursue for the main part of Darroch's rule, thus showing that $|\mu(\lambda,\bp)-m_\pm(\lambda,\bp)|\leq 1$ 
holds
in consequence of the chain of inequalities
$$
\ldots<\min_{{\cal F}_{k}}\mu(\lambda,\bp)=k-1+\frac{1}{k+1}<  k=  \max_{{\cal F}_{k}}\mu(\lambda,\bp)<\ldots
$$
In particular,
we obtain $\max_{{\cal F}_k}\mu(\lambda,\bp)=k$
by sending $n\to0$  in  (\ref{Dar}) and (\ref{maxmin}).
But this  implies existence of  bifurcation points on ${\cal M}_k$ that 
need to be explicitly found. Indeed, from the Poisson$(k)$ case we know that $(k,{\boldsymbol 0})\in {\cal M}_k \cap {\cal F}_k$ is one such point.

To complete the picture
we procede  on shoulders of Samuels \cite{Samuels} by generalising his Theorem 1.

\begin{theorem}\label{meanmode}
 If 
$k<\mu(\lambda,\bp)<k+1$ for some integer $k\geq0$
then 
$$
m_+(\lambda,\bp\setminus\max p_i)\leq k\leq m_-(\lambda,\bp)\leq m_+(\lambda,\bp)\leq k+1.
$$
If   $\mu(\lambda,\bp)=k$ then also $m_+(\lambda,\bp)=m_-(\lambda,\bp)=k$, unless
$(\lambda,\bp)=(k-a,\underbrace{1,\ldots,1}_{a},0,0,\ldots)$   {\rm (}for some $0\leq a<k, ~k\geq 1${\rm )}
in which case  the distribution is 
shifted Poisson  on $\{a,a+1,\ldots\}$
with the  twin mode
 $\{k-1,k\}$.
\end{theorem}

\proof The case $k=0$ being trivial, let $k\geq1$, and suppose in the first instance that $p_1<1$ (hence $p_i<1$ for all $i$).
Suppose $\sum_i p_i +\lambda>k$. To keep the balance in the identity of Lemma \ref{L1} (vi) we should have either 
\begin{equation}\label{Sam1}
f(k-1;\lambda,\bp)<f(k;\lambda,\bp)
\end{equation}
or $f(k-1;\lambda,\bp\setminus p_i)<f(k;\lambda,\bp)$ for some $i$. The second option is necessary if $\lambda=0$, and if it holds the recursion (i) of Lemma \ref{L1} implies
$f(k-1;\lambda,\bp\setminus p_i)<f(k;\lambda,\bp\setminus p_i)$, 
which is equivalent to $k\leq m_-(\lambda,\bp\setminus p_i)$. In the latter case the monotone likelihood property in Lemma \ref{L1} (vi) still leads to (\ref{Sam1}), whence
$k\leq m_-(\lambda,\bp)$ holds in any case.

Suppose $\sum_i p_i +\lambda<k+1$. By the same balance identity either $f(k;\lambda,\bp)>f(k+1;\lambda,\bp)$ (then $m_+(\lambda,\bp)\leq k$) or 
$f(k;\lambda,\bp\setminus p_i)>f(k+1;\lambda,\bp)$ for some $i$. In the second case we find as above that
$f(k;\lambda,\bp\setminus p_i)>f(k+1;\lambda,\bp\setminus p_i),$
whence by the monotone likelihood applied in the opposite direction
$f(k;\lambda,\bp\setminus \max p_i)>f(k+1;\lambda,\bp\setminus \max p_i).$
Thus in any case $ m_+(\lambda,\bp\setminus \max p_i)\leq k$, which in turn yields   $m_+(\lambda,\bp)\leq k+1$.

Now consider the case  $\sum_i p_i +\lambda=k\geq1$ and $\bp\neq0$. 
Suppose first that $f(k;\lambda,\bp)>f(k+1;\lambda,\bp)$, then $k\leq m_-(\lambda,\bp)$ and there exists $i$ such that
$f(k;\lambda,\bp)\leq f(k-1;\lambda,\bp\setminus p_i)$. 
Therefore $f(k;\lambda,\bp\setminus p_i)\leq f(k-1;\lambda,\bp\setminus p_i)$, which gives $k-1\geq m_-(\lambda,\bp\setminus p_i)$. Hence 
$k\geq m_+(\lambda,\bp)$ and putting things together   $m_+(\lambda,\bp)=m_-(\lambda,\bp)=k$.

Suppose $f(k;\lambda,\bp)\leq f(k-1;\lambda,\bp)$. Regardless of $\lambda$ there exists $i$ such that $f(k;\lambda,\bp)\geq f(k-1;\lambda,\bp\setminus p_i)$.
But then $f(k;\lambda,\bp\setminus p_i)\geq f(k-1;\lambda,\bp\setminus p_i)$ and by the monotone likelihood $f(k;\lambda,\bp)> f(k-1;\lambda,\bp)$,
which contradicts the assumption. Thus this case is excluded.

The case  $\bp=(k;{\boldsymbol 0})$ is  purely Poisson, with the mode $\{ k-1,k \}$.

Finally, we lift the restriction on parameters 
to allow some number $a$ of  $1$'s  among the $p_i$'s, $1\leq a<k$. The above conclusions are readily adjusted by 
reducing to the shifted  probability function $f(k-a;\lambda,\bp)$.
\endpf
\noindent
In the finitary case the first part of the argument applied to  the profile $\bp\cup \min p_i$ yields a slightly stronger implication
$$
\mu(\lambda,\bp)+\min p_i>k\Longrightarrow m_-(\lambda,\bp)\geq k,
$$
which is similar to   Equation (15) of  \cite{Samuels}.


For $k\geq 0$ let
$${\boldsymbol[}{\cal F}_{k},{\cal F}_{k+1}{\boldsymbol]}:=\{(\lambda,\bp)\in \nabla: f(k-1;\lambda,\bp)\leq f(k;\lambda,\bp)\geq f(k+1;\lambda,\bp)\},$$
which is the `interval' of parameters satisfying $k\in m(\lambda,\bp)$.

\begin{corollary} The only common points of some of the sets 
$ {\cal F}_{0}, {\cal M}_{0},  {\cal F}_{1}, {\cal M}_{1}, {\cal F}_{2},{\cal M}_{2},\ldots$
are the points of ${\cal M}_k \cap {\cal F}_k$ which correspond to shifted Poisson distributions with mean $k$. 
Every domain 
${\boldsymbol[}{\cal F}_{k},{\cal F}_{k+1}{\boldsymbol]}$
is  a connected compact,  which is  further divided in two connected components by  the convex compact ${\cal M}_k$.
\end{corollary}
\proof 
Every point $(\lambda,\bp)$ in the said domain can be connected to ${\cal M}_k$ by a segment of the ray through the point.
In more detail, 
if $(\lambda,\bp)$ satisfies  $f(k-1;\lambda, \bp)<f(k;\lambda, \bp)$ and $\mu(\lambda,\bp)<k$, then there exists $z>1$ such that the first inequality holds for $(z\lambda,z\bp)$
and $\mu(z\lambda,z\bp)=k$. If $\mu(\lambda,\bp)<k$ we choose deformation with $0<z<1$.
Furthermore,   any two points in  ${\cal M}_k$    can be connected by a linear segment because the set is convex.
Thus any two points in the domain can be connected by  a path with three contiguous segments.

The set  ${\boldsymbol[}{\cal F}_{k},{\cal F}_{k+1}{\boldsymbol]}$ is compact as being a closed subset of a part of $\nabla$ defined by a constraint of the form $a\leq \lambda\leq b$.
\endpf

A brief remark concerning the median.
Samuels and Jogdeo \cite{Jogdeo} showed for the finitary case that if the mean is $k$ then the median is also $k$.
This generalises by approximation to extended Bernoulli sums, since
the median $k^*$ is defined by the condition that each of the events $S\leq k^*$ and $S\geq k^*$ have probability at least $1/2$.

\section{Cross modal families}\label{S6}

With the above splitting of $\nabla$ in `intervals',
the modal ridge for the family of all extended binomial sums is
$$M=\bigcup\limits_{k=0}^\infty \{k\}\times {\boldsymbol[}{\cal F}_{k},{\cal F}_{k+1} {\boldsymbol]}.$$
A family of extended Bernoulli sums corresponds to some set of parameters ${\cal P}\subset\nabla$,
whose modal ridge is just the  slice-wise intersection
$$\bigcup\limits_{k=0}^\infty \{k\}\times    ({\cal P}\cap  {\boldsymbol[}{\cal F}_{k},{\cal F}_{k+1} {\boldsymbol]}).$$

For ease of exposition we shall consider families with the infinite range $R=\{0,1,\ldots\}$, leaving to the reader the case of a finite integer interval.
With this convention in mind, 
the family $\cal P$ is cross modal if for every $k\geq 0$
the maximum of likelihood $f(k;\cdot)$ in $\cal P$ is attained within ${\cal P} \cap{\boldsymbol[}{\cal F}_{k},{\cal F}_{k+1} {\boldsymbol]}$.
As an immediate consequence of the extended Darroch's rule in Theorem \ref{meanmode} we have
a result concerning the mean-matching families.

\begin{prop} \label{prop2}
If for every $k\geq0$ the maximum  $\max\{f(k;\lambda,\bp): (\lambda,\bp)\in{\cal P}\}$  is attained on  ${\cal P}\cap {\cal M}_k$ then the family $\cal P$ is cross modal.
\end{prop}

In more detail, the cross modality imposes the following  three conditions on $\cal P$:
\begin{itemize}

\item[(i)] the family is {\it unperforated}, that is ${\cal P}\cap  {\boldsymbol[}{\cal F}_{k},{\cal F}_{k+1} {\boldsymbol]}\neq\varnothing$,
\item[(ii)]  in ${\cal P}\cap{\boldsymbol[}{\cal F}_{k},{\cal F}_{k+1}{\boldsymbol]}$ the maximum of $h(\lambda,\bp)=f(k;\lambda,\bp)$
 is attained  at some point,
\item[(iii)] the maximum peak height in (ii) is the overall  largest value of $f(k;\lambda,\bp)$ for $(\lambda,\bp)\in {\cal P}$.
\end{itemize}
Note that by this definition the sequence of Poisson distributions with $\lambda=0,2,4,\ldots$ is unperforated,
though odd integers do not appear in the role of a leading mode.

Conditions (i) and (ii) are self-evident and will be 
 implied by other properties of $\cal P$
or taken in the sequel for granted. 
  In particular, an obvious sufficient condition for (ii) is that ${\cal P}$ is closed.

\begin{figure}
\centering
\includegraphics[scale=0.5]{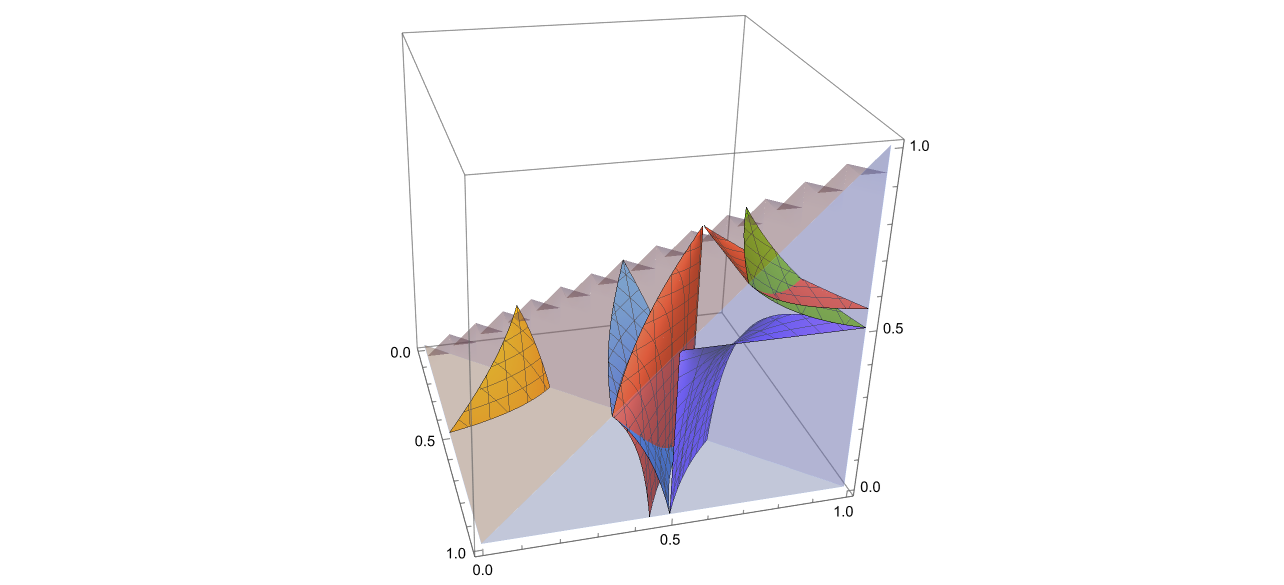}
\caption{For $n=3$, contour $f(2)=4/9$ (red colour) comprised of two branches meeting at $p_1=p_2=p_2=2/3$ and crossing both ${\cal F}_{2,3}$ and ${\cal F}_{3,3}$,
and contour $f(2)=1/2$ inside ${\boldsymbol[} {\cal F}_{2,3}, {\cal F}_{3,3}{\boldsymbol]}$.}
\label{CROSS1}
\end{figure}

Condition (iii) requires  much care, because  nonmodal likelihood under some distribution may be larger than the peak height of distribution with a larger mode,
for instance $f(2;1.6,{\boldsymbol 0})>f(2;2.5,{\boldsymbol 0})$.
A major complication for continuous-parameter families 
comes from the fact that 
 ${\cal F}_{k}$ is not a contour of the likelihood function. 
To the opposite, the contours of $f(k;\cdot)$
may cross  ${\cal F}_{k}$, or have branches separated by ${\cal F}_k$ or  lie completely  inside  ${\boldsymbol[}{\cal F}_{k},{\cal F}_{k+1}{\boldsymbol]}$,
see Figure \ref{CROSS1}. 
Even more weird,  in the $n=3$ finitary case there are  curves
 traversing ${\boldsymbol[}{\cal F}_{2,3},{\cal F}_{3,3}{\boldsymbol]}$ with increasing $\mu$
and $\max f(2;\cdot)$ lower than the  maximum of this function in neighbouring `intervals'.

The contours $f(k;\lambda,\bp)=z$ with $z\geq 1/2$ are enclosed in ${\boldsymbol[}{\cal F}_{k},{\cal F}_{k+1}{\boldsymbol]}$, thus 
offering an elementary sufficient condition for the case of  `big' peak heights.

\begin{prop} If  a closed family $\cal P$ 
for every $k\geq 0$ satisfies
$$\max\{ h(\lambda,\bp): (\lambda,\bp)\in{{\cal P}\cap {\boldsymbol[}{\cal F}_{k},{\cal F}_{k+1}{\boldsymbol]}}\}\geq \frac{1}{2}$$
then the family is cross modal.
\end{prop}
\proof The inequality
$f(k;\lambda,p)\geq 1/2$ implies that  $k$ is a mode, and that the likelihood maximum is attained in ${\cal P}\cap {\boldsymbol[}{\cal F}_{k},{\cal F}_{k+1}{\boldsymbol]}$.
\endpf

We take advantage of the natural order structure on $\nabla$ to introduce a reasonably large class of unperforated families ${\cal P}\subset\nabla$.
Define   ${\cal P}$ to be {\it directed} if $(0,{\boldsymbol 0})\in {\cal P}$ and the following conditions hold for $k\geq0$:

\begin{itemize}
\item[\rm (a)] for $(\lambda,\bp)\in {\cal P}$ with $m_+(\lambda,\bp)=k$ there exists $(\lambda',\bp')\in  {\cal P}$ such that 
$$(\lambda',\bp')\geq (\lambda,\bp),~ m_+(\lambda',\bp')=k+1~~{\rm and~~}f(k+1;\lambda',\bp')>f(k+1;\lambda,\bp),$$

\item[\rm (b)] for $(\lambda,\bp)\in {\cal P}$ with $m_-(\lambda,\bp)=k+1$ there exists $(\lambda',\bp')\in \cap {\cal P}$ such that 
$$(\lambda,\bp)\geq (\lambda',\bp'),~ m_-(\lambda',\bp')=k~~{\rm and~~}f(k;\lambda,\bp)<f(k;\lambda',\bp').$$

\end{itemize}
\noindent
Condition $(0,{\boldsymbol 0})\in {\cal P}$ can always be fulfilled by including zero in the family if needed.
The rationale behind the concept of directed family relates to the next property.
\begin{lemma} \label{L2} 
If $(\lambda',\bp')\geq (\lambda,\bp)$ for two distinct elements of $\nabla$, then 
 $k>m_+(\lambda',\bp')$ implies
$f(k;\lambda',\bp')>f(k;\lambda,\bp)$ provided $\lambda'+p_k'>0$.
\end{lemma}
\proof The last condition is needed to ensure $f(k;\lambda',\bp')>0$.
We may pass from $(\lambda,\bp)$ to  $(\lambda',\bp')$  
by successively increasing one of the components.
In the course of such possibly  infinitely many transitions the leading mode of intermediate vectors stays less than $k$, 
therefore  Lemma \ref{L1} (iii), (iv) ensures that the likelihood function $f(k;\cdot)$ strictly increases on some steps, whence eventually 
$f(k;\lambda',\bp')>f(k;\lambda,\bp)$. 
\endpf

\begin{theorem}\label{T3} Every closed directed family ${\cal P}\subset \nabla$ is cross modal.
\end{theorem}
\proof Condition $(0,{\boldsymbol 0})\in {\cal P}$ provides the base for induction to show that directed $\cal P$ is unperforated
in a strong sense, meaning that every possible value of $m_+$ and $m_-$ is taken.
Trivially $f(0;0,{\boldsymbol 0})=1$ is the absolute maximum of $f(0;\cdot)$. 

For $k>\ell\geq0$ let $(\lambda,\bp)\in {\cal P}$
have $m_+(\lambda,\bp)=\ell$. By induction and the assumption there exists a increasing
chain in $\cal P$ of distinct elements
$$(\lambda,\bp)\leq \ldots \leq (\lambda',\bp')\leq  (\lambda'',\bp''),$$
where the last two terms are chosen so that 
  $m_+(\lambda',\bp')=k-1, m_+(\lambda'',\bp'')=k$   and  $f(k;\lambda',\bp')< f(k;\lambda'',\bp'')$.
The likelihood function  $f(k;\cdot)$ is increasing strictly along the chain, as
the increase holds within the last two terms by the construction and 
for the remaining part of the chain  by Lemma \ref{L2}.

For $\ell>k$ a similar chain along which the likelihood decreases is  
constructed backwards from $(\lambda,\bp)\in {\cal P}\cap {\boldsymbol[}{\cal F}_{\ell},{\cal F}_{\ell+1}{\boldsymbol]}$.
It follows that the absolute maximum of $f(k;\cdot)$ can only  be found in ${{\cal P}\cap {\boldsymbol[}{\cal F}_{k},{\cal F}_{k+1}{\boldsymbol]}}\}$ and it exists 
since $\cal P$ is  closed.
\endpf

\subsection{Sequences}

The work of previous authors  \cite{Bruss,Cramer,  Moreno1, Moreno2, Tamaki} observed the cross-modality of the Bernoulli sums $S_n=B_1+\ldots+B_n$
for    $B_1, B_2,\ldots$ with given success probabilities, where  $n$ was treated as a parameter.
Bonferroni \cite{Bonferroni} was apparently the first to provide extended analysis of the peak height for the Bernoulli sums.

To embed this class of distributions in the general framework, one may consider a sequence
${\cal P}=\{(0,{\boldsymbol 0}),(\lambda_1,\bp_1), (\lambda_2,\bp_2),\ldots\}\subset \nabla$ such that
$(\lambda_{n+1},\bp_{n+1})$ is derived from its predecessor  $(\lambda_n,\bp_{n})$  by  one of the following operations:
\begin{itemize}
\item[{\rm (o1)}]  replacing  one of the success probabilities $p_i$ by a larger value $p_i'>p_i$, 
\item[{\rm (o2)}]  extending the vector by including a new component $p'$, 
\item[{\rm (o3)}]  incrementing the Poisson parameter by a quantity not exceeding $\gamma^*(\lambda_n,\bp_n)$,
\end{itemize}
where the parameters  of extension $i, p_i'$ or $p'$ may depend on $n$.
Obviously, (o1), (o3) correspond to convolution with a Bernoulli or Poisson variable, respectively.
The operation (o2) is a special case of (o1) `replacing zero by a Bernoulli variable'. The constraint in (o3) is needed to 
avoid decreasing likelihood.
The increments of the mean satisfy $\mu(\lambda,\bp_{n+1})-\mu(\lambda,\bp_{n})\leq 1$, whence by
Theorem \ref{meanmode}   ${\cal P}$ is unperforated. Applying Lemma \ref{L1} (iii) we see that ${\cal P}$ is directed, hence cross modal.
If  $\mu(\lambda,\bp_n)\to\infty$ as $n\to\infty$, the family has infinite range.

A well-known example is the Karamata-Stirling distribution, which we record with two parameters
\begin{equation}\label{KS}
f(k;n,t)= \left[  {n\atop k}  \right]\frac{t^k}{(t)_n}, ~~~t>0,~ 1\leq k \leq n.
\end{equation}
This is a Poisson-binomial distribution with the success probabilities $p_i=t/(t+i-1), i\geq 1$ (note: $p_1=1$).
The $t=1$ case relates to the random records model, and the general $t>0$ case plays a major role in the 
study of permutations, partitions and other combinatorial structures \cite{ABT}. 
With the aid of asymptotic expansion of the distribution, Kabluchko et al \cite{Kabluchko} 
have found approximation for the peak height and
showed that for large enough $n$ the mode is unique
and coincides with 
$   \lfloor u(n,t)-1/2\rfloor$ or     $ \lceil u(n,t)-1/2 \rceil\}$, where $u(n,t):=t\,(\log n + \partial_t \log \Gamma(t)),$
with the exact formula $m_{\pm}(n,t)= \lfloor u(n,t)\rfloor$ being valid for a set of $n$'s of asymptotic density $1$ (though not for all $n$).
As has been noticed in \cite{Pitman},
Darroch's rule applied to this instance 
gives  bounds for all $n$
$$    \lfloor \mu(n,t)\rfloor \leq m_{\pm}(n,t)\leq     \lceil \mu(n,t)\rceil ~~~{\rm for}~~\mu(n,t)=\sum_{i=1}^n \frac{t}{t+i-1},$$
and essentially reduces  the localisation to the asymptotic expansion of the shifted  harmonic numbers, for which we refer to $\mu(n,t)/t$ \cite{Flajolet}.
The  uniqueness of the mode in case of irrational $t$ is also immediate.
A new twist we wish to add is that  the maxima of (\ref{KS}) in $n$ are those values where the mode increases, 
that is the likelihood maximiser is the inverse function to the mode, 
like in the binomial scheme in Figure  \ref{RIDGE}.
We leave for now the interpretation through waiting times open.

There are many examples of row-wise totally positive 
triangular arrays of combinatorial numbers \cite{Pitman, Pollard}, 
leading to Bernoulli sums with parameters $\bp_n$  not explicitly known and not following the simple patterns (o1) or (o2). 
Stochastic monotonicity of  the associated sums $S_n$ is a necessary but not sufficient condition for Theorem \ref{T3} to work.
A numerical evidence  suggests that the array of normalised Stirling numbers of the second kind is cross modal.
See Rukhin et al \cite{Rukhin} for a related study of  the likelihood ratios for a  convolution of two binomial distributions.

\subsection{Power series and other continuous-parameter families}

For $a_0>0, r_i>0$ with $\sum_i r_i<\infty$ consider the entire function of the form
\begin{equation}\label{TPPS}
F(t) = \sum_{k=0}^\infty a_k t^k=a_0 e^{\lambda t}\prod_{i=1}^\infty (1+r_i t),
\end{equation}
which defines a power series family of logconcave distributions with parameter $t\geq0$.

The distribution with parameter $t=1$ is the law of an extended Bernoulli sum $S(\lambda,\bp)$ with Poisson parameter $\lambda$ and success probabilities $p_i=r_i/(1+r_i)$, as in (\ref{entireGF}).
For the generic $t\geq0$ the distribution  has probability generating function 
$$
\frac{F(tz)} {F(t)}=\frac{ a_0}{F(t)} e^{\lambda tz}\prod_{i=1}^\infty (1+r_i tz),
$$ 
corresponding to the extended Bernoulli sum with Poisson parameter $t\lambda$ and odds $tr_i$. 
This  family of power series distributions is cross modal by Proposition \ref{prop2} which in turn followed 
from the extended Darroch's rule in   
Theorem \ref{meanmode}. The theorem also implies that 
the criterion in Proposition \ref{prop1}holds
with  strict inequalities $t_k<\ell(k)<t_{k+1}$ and  $k-1<\mu(t_k)<k$ (last for $k\geq1$) 
Compare with  $\mu(t_k)=k$ in the pure  Poisson case.

Taking a geometric view will lead to a surprising conclusion, essentially matching
 the initial intuition that  `the mean and the mode grow at about the same pace'
was correct.
Note that  in the parametrisation by odds
the power series family corresponds to a ray.
In terms of the success probabilities the family has parameters $(\lambda(t), \bp(t))$, where
\begin{equation}\label{exten-par}
 \lambda(t)=\lambda t,~~{\rm and}~~p_i(t)=\frac{tp_i}{1-p_i+tp_i}
\end{equation}
are increasing functions.
Geometrically, $(\lambda(t), \bp(t)), t\geq0,$ is a continuous curve in $\nabla$, which is increasing (in the partial order) and unperforated since $\mu(\lambda(t),\bp(t))\uparrow\infty$.
The family is directed. Indeed, interpret  $t_k$ and $t_{k+1}$ as the `time'  when the curve enters into, respectively, exits from
${\boldsymbol[}{\cal F}_{k},{\cal F}_{k+1}{\boldsymbol]}$.
Lemma \ref{L2} tells us that $f(k;\lambda(t), \bp(t))$ strictly increases for $t\leq t_k$ and strictly decreases for $t\geq t_{k+1}$
thus the likelihood attains its absolute maximum inside ${\boldsymbol[}{\cal F}_{\ell},{\cal F}_{\ell+1}{\boldsymbol]}$.
Moreover by Lemma \ref{L1} (iii), we have
 $\partial_t f(k;\lambda(t_k), \bp(t_k))>0$ and $\partial_t f(k;\lambda(t_{k+1}), \bp(t_{k+1}))<0$, understood as transversal crossing the contours
of the likelihood function.
It again follows that the family of extended Bernloulli sums defined by the totally positive power series (\ref{TPPS}) is cross modal.
But neither general criterion for logconcave families in Proposition \ref{prop1} (iii) nor the proof of Theorem \ref{T3}   rely on the Darroch's rule, therefore the rule becomes a consequence. So to summarise:
\vskip0.3cm
\noindent
{\it since every extended Bernoulli sum is a member of a power series family {\rm (\ref{exten-par})}, the cross modality of the family implies the extended Darroch's rule}.
\vskip0.3cm
\noindent
Scrutinising this way of reason, three ingredients have been used: matching the mean 
by the likelihood maximisation 
for  power series families, including distribution in such a family and the elemenrary 
recursion in Lemma \ref{L1} (i) observed long ago by the  pioneers \cite{Bonferroni, Chebyshev}.

The roots of  some totally  positive entire functions (\ref{TPPS})  are known explicitly. 
For the example in Figure \ref{COSH} the Weierstrass product formula applied to  $\cosh{\sqrt t}$ gives the parameters of extended Bernoulli sum as $\lambda=0$ and
$$p_i=\frac{4}{4+(1-2i)^2\pi^2}, ~~i\geq 1.$$
Karamata-Stirling family of distributions (\ref{KS}) relates to a totally positive polynomial, for which we conclude on the cross modality  now in the variable $t$.

\section{Breaking the mode}\label{S7}

The value of Darroch's rule for localising the mode lies in the simplicity of the mean as 
compared to the whole multiparameter profile $(\lambda,\bp)$, 
which may not be known explicitly.
Bur even with  $(\lambda,\bp)$ given, it looks hard to control the mode under the (o2)-deformations of the profile.
From Theorem \ref{T1} we only  infer  that increasing  a particular $p_i$ will break the mode  if the success probability is  smaller than a certain bound 
depending on $(\lambda,\bp)$.

The peak skewness is an alternative characteristic of stability of the mode. 
Recall that  $\gamma^*(\lambda,\bp)$ is the minimal parameter of Bernoulli variable  needed to transport the mode to larger position.
 The monotonicity properties of $\gamma^*(\lambda,\bp)$ agree with its interpretation as a size of effort to break the mode, that is 
$$\partial_{p_i} \gamma^*(\lambda,\bp)<0, ~~~\partial_{\lambda} \gamma^*(\lambda,\bp)<0$$
in consequence of monotonicity of the likelihood ratio and the formula
$$\gamma^*(\lambda,\bp)= \left( 1+   \frac{1-\frac{f(m-1)}{f(m)} }{ 1-\frac{f(m+1)}{f(m)}           }          \right)^{-1}, $$
where for shorthand
$m=m_+(\lambda,\bp)$.
Curiosly, there is some connection with the mean.
In particular, if $(\lambda,\bp)\in {\cal M}_k$ one needs to increase the mean by at least $1/(k+2)$ to increment the leading mode by one,
hence $\gamma^*(\lambda,\bp)\geq 1/(k+2)$.
On the other hand, if $(\lambda,\bp)\in {\cal F}_k$ (the twin mode case),  to  increment the leading mode the mean should be increased by 
at least $1-1/(k+1)$ while $\gamma^*(\lambda,\bp)=1$.

\subsection{Two Bernoulli increments versus one}

The above  begs a question: how to transport the mode of $S$ by adding independent Bernoulli or Poisson variables  in the most economic way, as measured by the mean?

To that end,
suppose $S$ has probability function $f(k)$ with leading mode $m$, and we add an independent variable $Z$. A minute thought shows that shifting the mode to $m+1$ amounts to a linear programming problem 
with the balance equation as a constraint
$$
{\mathbb P}[S+Z=m]={\mathbb P}[S+Z=m+1],
$$
and the objective functional ${\mathbb E}Z$, where $Z$ can be restricted to have the support within $\{0,1,\ldots,m+1\}$.
If we only consider $Z$ with a
two-point distribution ${\mathbb P}[Z=s]=\delta$, ${\mathbb P}[Z=0]=1-\delta$ with given $s$, then 
the minimal necessary value of $\delta$  is found from the balance equation  as 
$$
\delta=\frac{\Delta f(m)}{ \Delta f(m) -\Delta f(m-s)},
$$
(generalising $\gamma^*(\lambda,\bp)$ in the $s=1$ case) which incurs the transportation cost $s\delta$.
For convexity reasons, if no  constraint on the distribution of $Z$ is imposed, the optimal distribution is two-point, with the range parameter found as the minimiser in 
\begin{equation}\label{min-s}
\min_s \frac{s \Delta f(m)}{ \Delta f(m) -\Delta f(m-s)},
\end{equation}
and $\delta=\delta(s)$ determined as above.

Within the class of extended Bernoulli sums, for the sake of optimisation it is sufficient to consider $Z$ representable as a Bernoulli sum with at most $m+1$ terms.
In the case $m=0$ the optimal transport involves a Bernoulli variable $Z$ with the success probability $\gamma^*(\lambda,\bp)$.
For $m>0$ we are willing to show that this one-term solution might not be optimal. 
Indeed, comparing $s=1$ and $s=2$ instances of (\ref{min-s}), the second will be more less costly if
$$ \frac{2 \Delta f(m)}{ \Delta f(m) -\Delta f(m-2)}<\frac{\Delta f(m)}{ \Delta f(m) -\Delta f(m-1)},$$
which is equivalent to
\begin{equation}\label{delta-3}
A:=\Delta^3 f(m-2)<0.
\end{equation}
 Condition  (\ref{delta-3}) is necessary for a two-term Bernoulli sum with some success probabilities $\alpha_1,\alpha_2$ to outperform a single Bernoulli with parameter $\gamma^*(\lambda,\bp)$.
Lemma \ref{L1} (viii) entails the expansion
\begin{equation}\label{dif-exp}
f(k;\lambda,\bp\cup\alpha_1,\alpha_2)=f(k;\lambda,\bp)-(\alpha_1+\alpha_2) \Delta f(k;\lambda,\bp)+\alpha_1 \alpha_2 \Delta^2 f(k;\lambda,\bp),
\end{equation}
which allows us to write the balance equation
$$f(m;\lambda,\bp\cup\alpha_1,\alpha_2)=f(m+1;\lambda,\bp\cup\alpha_1,\alpha_2)$$
in the form
\begin{equation}\label{quad}
C+ (\alpha_1+\alpha_2) B+ \alpha_1\alpha_2 A=0, ~~~{\rm where}~~C:= \Delta f(m), B:= -\Delta^2 f(m-1),
\end{equation}
and we have $C>0, B<0$ since $m$ is the mode.
In this notation, the peak skewness becomes $\gamma^*(\lambda,\bp)=-C/B$, which substituted in the balance equation yields 
$$
C+(\alpha_1+\alpha_2)B+ \alpha_1\alpha_2 C= C+\gamma^*(\lambda,\bp) B,
$$
and therefore
$$\alpha_1+\alpha_2 -\gamma^*(\lambda,\bp)=-\alpha_1\alpha_2\,\frac{A}{B}<0,$$
showing that two Bernoulli variables may be better than one.
To find the desired parameters explicitly,  set $\alpha:=\alpha_1=\alpha_2$ and solve the balance  equation $C+2B\alpha+A\alpha^2=0$ as
$$
\alpha= \frac{-B-\sqrt{B^2-AC}}{A}. 
$$
It is readily checked that  $2\alpha<-C/B$ indeed holds and thus $0<\alpha<1$.
In accord with \cite{Chebyshev} the optimum has been found in the class of  binomial distributions.

Thus condition (\ref{delta-3}) turns also sufficient. We have found numerically examples of  $S$ with $A<0$, for instance the Binomial$(12, 0.385)$ distribution has $m=5$ and $A=-0.003$,
which gives a tiny advantage to perturbation of $Z$ by a single term.

The form of balance equations with $Z$ being a sum of  $d\geq 3$ independent Bernoulli variables can be guessed from (\ref{dif-exp}) and (\ref{quad}).
However, the case of two terms suggests that  to get some further improvement  $\Delta^d f(m-d+1)$ should disagree with the sign-alternating 
pattern  emerging by the normal approximation.

\subsection{Breaking the Poisson mode}

Consider $S\stackrel{d}={\rm Poisson}(t), ~~m<t<m+1$. Adding $Z\stackrel{d}{=}$Poisson$(k+1-t)$  shifts the mode to $m+1$,
but adding a Bernoulli variable with success probability 
$$
\gamma^*(t,{\boldsymbol 0})= \frac{(m+1)t-t^2}{2(m+1)t-t^2-m(m+1)}
$$
offers a more spare mode transport with   $\gamma^*(t,\bp)<m+1-t$, the inequality holding 
since the cubic polynomial 
$$((m+1)t-t^2 )- [(2(m+1)t-t^2-m(m+1))(m+1-t)]$$
is strictly positive in the interval $(m, m+1)$.
From 
$\Delta^3 f(m-2)>0$ follows that no two-term deformation transporting the mode could have smaller mean.

\vskip0.5cm
\noindent
{\bf Acknowledgement} The author is indebted to Alex Marynych, Jim Pitman and Michael Farber for inspiring discussions, help with terminology and pointers to the literature.

\end{document}